\documentclass[12pt,reqno]{amsart}
\usepackage{amssymb,delarray}
\usepackage{amsfonts,amscd}
\usepackage{epsfig}
\usepackage[all]{xy}

\def\leq{\leqslant}
\def\geq{\geqslant}

\newtheorem{thm}{Theorem}
\newtheorem{lem}
{Lemma}
\newtheorem{prop}
{Proposition}
{Problem}
{Claim}
\newtheorem{df}
{Definition}
{Corollary}
\newtheorem{rem}
{Remark}
{Question}
\newtheorem{ex-thm}{Theorem-Example}

{\catcode`\@=11
\gdef\n@te#1#2{\leavevmode\vadjust{%
 {\setbox\z@\hbox to\z@{\strut#1}%
  \setbox\z@\hbox{\raise\dp\strutbox\box\z@}\ht\z@=\z@\dp\z@=\z@%
  #2\box\z@}}}
\gdef\leftnote#1{\n@te{\hss#1\quad}{}}
\gdef\rightnote#1{\n@te{\quad\kern-\leftskip#1\hss}{\moveright\hsize}}
\gdef\?{\FN@\qumark}
\gdef\qumark{\ifx\next"\DN@"##1"{\leftnote{\rm##1}}\else
 \DN@{\leftnote{\rm??}}\fi{\rm??}\next@}}

\begin{document}
\baselineskip=13.7pt plus 2pt 

\title[Rigid
germs of finite morphisms
]{Rigid
germs of finite morphisms of smooth surfaces and rational Belyi pairs
}

\author[Vik.S. Kulikov]{Vik.S. Kulikov}

\address{Steklov Mathematical Institute of Russian Academy of Sciences, Moscow, Russia}
 \email{kulikov@mi.ras.ru}

\dedicatory{} \subjclass{}
\thanks{ This work is supported by the Russian Science Foundation under grant no.  19-11-00237.}


\maketitle

\def\st{{\sf st}}

\quad \qquad \qquad

\begin{abstract} In \cite{K-rig}, a map $\beta:\mathcal R\to\mathcal{B}el$  from the set $\mathcal R$ of equivalence classes of rigid germs of finite morphisms  branched in germs of curves having  $ADE$ singularity types onto the set $\mathcal{B}el$  of rational Belyi pairs   $f:\mathbb P^1\to\mathbb P^1$ considered up to the action of $PGL(2,\mathbb C)$ was defined. In this article, the inverse images of this map are investigated  in terms of  monodromies of Belyi pairs.
\end{abstract} \vspace{0.5cm}

\def\st{{\sf st}}

\setcounter{section}{-1}

\section{Introduction}
In this article, we continue the investigate of properties of germs $F:(U,o')\to (V,o)$ of finite morphisms of smooth surfaces (further, for short, the {\it germs of covers}) started in \cite{K-4} and \cite{K-rig}. In \cite{K-4}, the notion of deformation equivalence of germs of covers was introduced.
A germ of cover $F:(U,o')\to (V,o)$ is called {\it rigid} if any  deformation equivalent to $F$ germ of cover $F_1:(U_1,o'_1)\to (V,o)$ is equivalent to it, that is, in short,  the covers $F$ and $F_1$ are the  covers different from each other on the change of coordinates in $(U,o')$ and $(V,o)$. In \cite{K-rig}, it was proved that
if the germ $(B,o)\subset (V,o)$ of the branch curve of a germ of cover $F:(U,o')\to (V,o)$ has one of $ADE$-singularity types, then  $F$ is a rigid germ.

Denote by $\mathcal R=(\bigcup_{n\geq 1}\mathcal R_{{\bf A}_n})\cup (\bigcup_{n\geq 4}\mathcal R_{{\bf D}_n})\cup(\bigcup_{n\in\{6,7,8\}}\mathcal R_{{\bf E}_n})$ the set of rigid germs of covers branched in curve germs having, resp., the singularity types ${\bf A}_n$, $n\geq 1$, ${\bf D}_n$, $n\geq 4$,  and ${\bf E}_6$, ${\bf E}_7$, ${\bf E}_8$.

A germ of cover $F$ of degree $\deg F=d$ defines a homomorphism $F_*:\pi_1(V\setminus B,p)\to \mathbb S_d$ (the {\it monodromy} of the germ $F$), where $\mathbb S_d$ is the symmetric group acting on the fibre $F^{-1}(p)$.  The group $G_F=\text{im} F_*\subset \mathbb S_d$ is called the ({\it local}) {\it monodromy group} of $F$. Note that $G_F$ is a transitive subgroup
of $\mathbb S_d$. By  Grauert - Remmert - Riemann - Stein Theorem (\cite{St}), the monodromy homomorphism $F_*$ defines a cover $F$ uniquely up  to equivalence.

Denote by $\mathcal{B}el$ the set of rational Belyi pairs considered up to the action of the group $\text{PGL}(2,\mathbb C)$ on $\mathbb P^1$. A cover $f:\mathbb P^1\to\mathbb P^1$, defined over the algebraic closure $\overline{\mathbb Q}$ of the field of rational numbers $\mathbb Q$,
is called a {\it Belyi pair} if it is branched at most in three points, $\mathcal Bel=\mathcal Bel_{2}\cup \mathcal Bel_3$, where  $\mathcal Bel_{2}$ is the set of Belyi pairs branched at most in two points and the Belyi pairs $f\in \mathcal Bel_3$  branched in three points. Further,
we will assume that $f\in \mathcal Bel_{2}$ a cover given in nonhomogeneous coordinates by function  $z=x^n$, $n\geq 1$, and  its branch locus  is $B_f=\{ 0,\infty\}$ (if $n\geq 2$), and  the branch locus of $f\in \mathcal Bel_{3}$ is
$B_f=\{ 0,1,\infty\}$.

In \cite{K-rig}, it was defined a map $\beta:\mathcal R\to\mathcal Bel$ as follows. Let $F:(U,o')\to (V,o)$ be a germ of  cover branched in a germ $(B,o)\subset (V,o)$ having one of $ADE$ singularity types and let $\sigma:\widetilde V\to V$ be a sequence of $\sigma$-processes with centers at points such that $\sigma^{-1}(B)$ is a divisor with normal crossings  (but, if the singularity type of $B$ is ${\bf A}_0$ or ${\bf A}_1$ then $\sigma$ is the single $\sigma$-process with center at $o$). Denote by $E\subset \widetilde V$ the exceptional curve of the last $\sigma$-process
and by $\widetilde F: \widetilde U\to \widetilde V$ and $\tau:\widetilde U\to U$ two natural holomorphic maps from the normalization  of the fibre product $\widetilde U=U\times_V \widetilde V$ of the holomorphic maps $F:(U,o')\to (V,o)$ and $\sigma:\widetilde V\to (V,o)$. It is easy to show that $C=\widetilde F^{-1}(E)$ is an irreducible rational curve and the restriction $f=\widetilde F_{\mid C}: C\to E$ is branched at most in three points. By definition, the map $\beta$ sends $F\in \mathcal R$ to $f\in \mathcal {B}el$,

\begin{picture}(300,95)
\put(25,85){$\phantom{aaa}$} \put(25,70){$\mathbb P^1\simeq C$}
\put(68,70){$\subset$}
\put(85,70){$\widetilde U$}
\put(105,75){\vector(1,0){60}}
\put(125,80){$\tau$}
\put(170,70){$U\ni o'$}
\put(25,-5){$\mathbb P^1\simeq E$}
\put(68,-5){$\subset$}
\put(85,-5){$\widetilde V$}
\put(107,-0){\vector(1,0){60}}
\put(125,-10){$\sigma$}
\put(170,-5){$V\ni o$}
\put(28,35){$f$}
\put(43,62){\vector(0,-1){50}}
\put(93,35){$\widetilde F$}
\put(90,62){\vector(0,-1){50}}
\put(180,35){$F$}
\put(175,62){\vector(0,-1){50}}
\put(250,35){$\beta: F\in\mathcal R \mapsto \beta(F)=f\in\mathcal{B}el$.}
\end{picture}\vspace{0.5cm}

Similar to the two-dimensional case, a cover $f\in\mathcal Bel$ defines a homomorphism $f_*:\pi_1(\mathbb P^1\setminus B_f,p) \to \mathbb S_n$ (the {\it monodromy} of $f$), where $n=\deg f$. The image $G_f=\text{im} f_*\subset \mathbb S_n$ is called the  {\it monodromy group} of $f$. If $f\in\mathcal Bel_{2}$, then $G_f= \mu_n\subset\mathbb S_n$ is a cyclic group of order $n$.

The group $\pi_1(\mathbb P^1\setminus \{ 0,1,\infty\},p)$ is the free group generated by two simple
loops $\gamma_0$ and $\gamma_1$ around the points $0$ and $1$ such that the loop $\gamma_{\infty}=\gamma_0\gamma_1$ is the trivial element in $\pi_1(\mathbb P^1\setminus \{ 0,1\},p)$.
For $f\in \mathcal Bel_3$ denote by $$T_c(f)=\{ c_i=(m_{1,i},\dots ,m_{k_i,i})\}_{m_{1,i}+\dots + m_{k_i,i}=\deg f,\,\, i\in \{ 0,1,\infty\}}$$
the set of cycle types of permutations $f_*(\gamma_i)$. 
Then, by Hurwitz formula connecting the degree of $f:\mathbb P^1\to \mathbb P^1$ and the orders of ramification at the critical points of $f$, we have the following equality:
\begin{equation} \label{H} n+2=k_0+k_1+k_{\infty}.\end{equation}
Conversely, if a transitive group $G\subset \mathbb S_n$ is generated by two permutations $\sigma_0$ and $\sigma_1$ such that their cycle types and the cycle type of $\sigma_{\infty}=\sigma_0\sigma_1$ satisfy equation (\ref{H}) then there is a rational Belyi pair $f$ such that $f_*(\gamma_i)=\sigma_i$.

In \cite{K-rig}, it was shown that for $F\in \mathcal R$ the covers $\widetilde F$ and $F$ can be represented as  compositions of two finite maps (see diagram $(*)$ in Subsection \ref{neces}), $\widetilde F=\widetilde H_2\circ\widetilde H_1$ and $F=H_2\circ H_1$, where $\widetilde H_1:\widetilde U\to \widetilde W$ and  $H_1:U\to W$ are cyclic covers (here $\widetilde W$ and $W$ are normal surfaces) such that $\widetilde H_{1\mid C}: C\to \widetilde H_1(C)$ is an isomorphism, and $\widetilde H_2:\widetilde W\to \widetilde V$ and $H_2:W\to V$ are finite covers such that the monodromy homomorphisms $\widetilde H_{2*}$ and $H_{2*}$ of the covers $\widetilde H_2$ and$H_2$) can be identified with the monodromy homomorphism $\beta(F)_*$ of the Belyi pair $\beta(F)$.

In  Section \ref{section2} (see Theorem \ref{main2}), a description of $\mathcal R_T\cap\beta^{-1}(f)$ is given for all $f\in\mathcal{B}el$ and $\mathcal R_T\subset \mathcal R$
in terms of the monodromy homomorphism $f_*$. In particular, in Section \ref{thmD_4}, we prove

\begin{thm}\label{B-D_4} Let $f\in \mathcal{B}el$, $\deg f=n>1$ and $B_f\subset \{0,1,\infty\}$, be given by two coprime homogeneous in variables $x_1,x_2$ forms $(h(x_1,x_2):h_2(x_2,x_2))$,
$$f:\,\, (x_1:x_2)\mapsto (h(x_1,x_2):h_2(x_2,x_2),$$
and let $p_1=(0,1)$ and $p_2=(1,0)$ be two points such that $\{ f(p_1), f(p_2)\}\cup B_f=\{ 0,1,\infty\}$. Then a cover
$F:(U,o')\to (V,o)$ given by functions
\begin{equation} \label{b-d_4}
u=h_1(z^{m_1},w^{m_2}),\quad v=h_2(z^{m_1},w^{m_2}) \end{equation}
belongs to $\mathcal R_{{\bf D}_4}$, where $m_1, m_2\in \mathbb N$ such that $G.C.D.(m_1,m_2)=1$ and where  $m_1>1$ and $f(p_1)=1$ if $f\in \mathcal{B}el_2$.

Conversely, any $F\in\mathcal R_{{\bf D}_4}$ is equivalent to a cover given by functions of the form {\rm (\ref{b-d_4})} and its image $\beta(F)$ is  $f:\,\, (x_1:x_2)\mapsto (h(x_1,x_2):h_2(x_2,x_2))$.
\end{thm}

A complete description of the sets $\mathcal R_T\cap\beta^{-1}(\mathcal{B}el_2)$ for $\mathcal R_T\subset \mathcal R$ is given in

\begin{thm} \label{beta2} If $F\in (\bigcup_{k=1}^{\infty}\mathcal R_{{\bf A}_{2k}})\cup\mathcal R_{{\bf E}_6}\cup\mathcal R_{{\bf E}_8}$, then $\beta(F)\in\mathcal{B}el_3$.

If $\beta(F)=f\in\mathcal{B}el_2$, $\deg f=n$, for $F\in \mathcal R\setminus((\bigcup_{k=1}^{\infty})\mathcal R_{{\bf A}_{2k}})\cup\mathcal R_{{\bf E}_6}\cup\mathcal R_{{\bf E}_8})$, then $F$ is equivalent to one of the following covers: \newline
$F\in \mathcal R_{{\bf A}_0}:\,\, \qquad \qquad \qquad u=z^m,\,\, v=w$,

where $m\geq 1$, $n= 1;$
\newline
 $F\in \mathcal R_{\bf A_1}$:\,\, $\qquad \qquad \qquad\, u=z^{nm_1},\,\, v=w^{nm_2}$,

where $n\geq 1, m_1\geq m_2\geq 1$; 
\newline
$F\in \mathcal R_{\bf A_{2k+1}}$, $k\geq 1$:\,\, $\qquad u=(z^m+w^{m_0})^n$, \,\, $v=w$,

where $n,m,m_0>1$, $k+1=nm_0$; 
\newline
$F\in\mathcal R_{\bf A_{2k+1}}$, $k\geq 1$: $\qquad u=z^{nm_1},\,\, v=z^{m_1}+w^{m_2}$,

where
$m_1\geq 1$, $n,m_2>1$; 
\newline
$F\in\mathcal R_{\bf A_{2k+1}}$, $k\geq 1$: $\qquad u=(\omega_jz^{m_1}-w^{m_2})^n,\,\, v=z^{m_1}-w^{m_2},$

where $n, m_1,m_2>1$, 
$\omega_j=exp({2\pi j i/n})$, $1\geq j\geq n-1$; \newline
$F\in\mathcal R_{\bf D_{2k+3}}$, $k\geq 1$: $\qquad u=z^{2m_1}, \,\, v=z^{m_1(2k+1)}+w^{m_2}$,

where $m_1\geq 1$, $m_2>1$, $n=2$,  
 $\text{\rm G.C.D.}(2k+1,m_2)=1$;  \newline
$F\in\mathcal R_{\bf D_{2k+2}}$, $k\geq 2$: $\qquad u=z^{n_1m_1},\quad v=(z^{m_1k_2}+w^{m_2})^{n},$ 

where $k=k_1k_2,\, n=n_1k_1\geq 2$, $m_1,\, m_2\geq 1$,  
$\text{\rm G.C.D.}(nm_2,k_2)=1$; \newline
$F\in\mathcal R_{\bf D_{2k+2}}$, $k\geq 2$: $\qquad u=(z^{m_1}-w^{m_2})^{n_1},\quad v=z^{m_1n},$ 

where $n=n_1k\geq 1$, $m_1,m_2\geq 1$; 
\newline
$F\in \mathcal R_{\bf D_{2k+2}}$, $k\geq 2$: $\qquad u=(z^{m_1}-w^{m_2})^{n_1},\,\, v=(z^{m_1}-\omega_jw^{m_2})^{n},$ 

where  $n=n_1k\geq 2$, $m_1,m_2\geq 1$, 
$\omega_j=exp(2\pi ji/n)$, $j=1,\dots, n-1$; \newline
$F\in\mathcal R_{\bf D_{4}}$: $\qquad \qquad \qquad u=z^{m_1n},\,\, v=(z^{m_1}+w^{m_2})^n,$

where $n\geq 2$, $m_1$, $m_2\geq 1$; 
\newline
$F\in\mathcal R_{\bf D_{4}}$: $\qquad \qquad \qquad u=(z^{m_1}-w^{m_2})^n,\,\, v=(z^{m_1}-\omega_jw^{m_2})^n,$ 

where $n\geq 2$, $m_1$, $m_2\geq 1$, 
$\omega_j=exp(2\pi ji/n)$, $1\leq j\leq n-1$; \newline
$F\in\mathcal R_{\bf E_{7}}$: $\qquad \qquad \qquad u=z^{3m_1},\,\, v=z^{2m_1}+w^{m_2},$ 

where $m_1\geq 1$, $m_2>1$.

In all cases $\text{\rm G.C.D.}(m_1,m_2)=1$.
\end{thm}
The proof of Theorem \ref{beta2} is given in Section \ref{Bel2}.

\section{Preliminary results} \label{preliminary}
\subsection{On the fundamental groups.} Denote by $(X,o)$  a germ of normal surface and $(B,o)=\bigcup_{j=1}^m B_j$ a union of $m\geq 0$ irreducible curve germs $(B_j,o)\subset (X,o)$.  Let  $\sigma:\widetilde X\to (X,o)$ be the {\it minimal resolution of the singularity} of the pair $(X,B,o)$, that is, $\widetilde X$ is smooth and $\widetilde B=\sigma^{-1}(B)$ is a divisor with normal crossings in which each $(-1)$-curve intersects at least with three irreducible components of $\sigma^{-1}(B)$. Below, we will assume that  $\sigma^{-1}(o)=\bigcup_{j=1}^kE_j$ is a union of rational curves and the dual graph of $\sigma^{-1}(o)$ is a tree. Also, if this does not lead to a misunderstanding, the proper inverse images $\sigma^{-1}(B_j)$ of the irreducible curve germs $B_j$ of $(B,o)$ will be marked with the same letter $B_j$.

The dual weighted  graph $\Gamma(\widetilde B)$  of  $\widetilde B$ is a tree having $m+k$ vertices $v_j$. The vertices $v_j$, $j=1,\dots, m$,  correspond to the curve germs $B_j$ and their weights are $w_j=0$, the vertices $v_{m+j}$, $j=1,\dots, k$,  correspond to the curves   $E_{j}$ and their weights are $w_{m+j}=-(E_{j}^2)_{\widetilde X}$.
For each pair of vertices $v_i$ and $v_j$ of $\Gamma(\widetilde B)$, we define
$$ \delta_{i,j}=\left\{ \begin{array}{ll} 1, & \text{if}\,\, v_i\,\, \text{and}\,\, v_j\,\, \text{are connected by an edge in}\,\, \Gamma(\widetilde B),
\\ 0, & \text{if}\,\, v_i\,\, \text{and}\,\, v_j\,\, \text{are not connected by an edge in}\,\, \Gamma(\widetilde B), \\ 0, & \text{if}\,\, i=j.\end{array} \right. $$

The following Theorem \ref{Mum} allows to obtain a presentation of the fundamental group $\pi_1(\widetilde X\setminus \widetilde B)$ in terms of the graph
$\Gamma(\widetilde B)$. The proof of this Theorem coincides almost word for word with the proof of similar statement in \cite{Mu} (see also \cite{K-rig}) and therefore it will be omitted.

\begin{thm}\label{Mum} 
The group $\pi_1(\widetilde X\setminus \widetilde B)$ is generated by $m+k$ elements $b_1,\dots, b_m$   being in one-to-one correspondence with the vertices $v_1,\dots, v_m$
and  $e_{m+1},\dots, e_{m+k}$  being in one-to-one correspondence with the vertices $v_{m+1},\dots ,v_{m+k}$ of $\Gamma(\widetilde B)$, and  being subject to the following defining  relations:
$$ \begin{array}{rll}
e_{m+i}^{-w_{m+i}}\cdot b_1^{\delta_{1,m+i}}\cdot .\, .\, .
\cdot b_m^{\delta_{m,m+i}}\cdot
e_{m+1}^{\delta_{m+i,m+1}}\cdot . \, .\, .\cdot
e_{m+k}^{\delta_{m+i,m+k}}
&  =1\quad & \text{for}\,\,i =1,\dots, k, \\
{[}b_j, e_{m+i}{]} & =1 \quad & \text{if}\,\,\,\, \delta_{j,m+i}=1,     \\
{[}e_{m+i_1}, e_{m+i_2}{]} & =1 \quad & \text{if}\,\,\,\, \delta_{m+i_1,m+i_2}=1.
\end{array}
$$
\end{thm}
\begin{rem} The generators  $b_1,\dots, b_m$ and  $e_{m+1},\dots, e_{m+k}$ of $\pi_1(\widetilde X\setminus \widetilde B)$ in Theorem {\rm \ref{Mum}} are presented by some simple loops around the curves corresponding to them {\rm (see \cite{K-rig})}. \end{rem}

The following Lemma is well-known (see, for example, \cite{K-Sh}).
\begin{lem} \label{sigma} Let $(Y,o)$ be a germ of smooth surface, $\sigma: X\to (Y,o)$ the $\sigma$-process with center at $o$, $(C_1,o)$ and $(C_2,o)$ two smooth curve germs in $(Y,o)$ meeting transversally at $o$. Then  $\gamma_E=\gamma_1\gamma_2$ in $\pi_1(Y\setminus (C_1\cup C_2))\simeq \pi_1(X\setminus \sigma^{-1}(C_1\cup C_2))$, where $\gamma_E$ is an element in $\pi_1(X\setminus \sigma^{-1}(C_1\cup C_2))$ represented by a simple loop around the exceptional curve $E=\sigma^{-1}(o)$ and $\gamma_j$, $j=1,2$, are the elements
represented by  simple loops around $C_j$.
\end{lem}
\subsection{Graphs of resolution of singularities of $ADE$ singularity types} \label{graphade}
Remind that the equations of curve germs $(B,o)$ having one of $ADE$ singularity types are the following (\cite{Ar}) :
\begin{itemize}
\item[\phantom{a}]  ${\bf A}_n: \qquad u^2-v^{n+1}=0,\,\, n\geq 0$;
\item[\phantom{a}]  ${\bf D}_n: \qquad v(u^2-v^{n-2})=0,\,\, n\geq 4$;
\item[\phantom{a}]  ${\bf E}_6: \qquad u^3-v^4=0$;
\item[\phantom{a}]  ${\bf E}_7: \qquad u(u^2-v^3)=0$;
\item[\phantom{a}]  ${\bf E}_8: \qquad u^3-v^5=0$.
\end{itemize}

The graph $\Gamma(\widetilde B)$ of the curve germ $(B,o)$ of singularity type ${\bf A}_{2k+1}$, $k\geq 0$, is depicted on Fig. 1 (if $k=0$ then the weight  of the vertex $e_3$ is equal to $-1$).

\begin{picture}(300,85)
\put(85,30){\circle*{3}}\put(80,35){$\mbox{}_{-2}$}
\put(82,19){$\mbox{e}_{3}$}\put(85,30){\line(1,0){40}}
\put(125,30){\circle*{3}}\put(120,35){$\mbox{}_{-2}$}
\put(122,19){$\mbox{e}_{4}$}\put(125,30){\line(1,0){40}}
\put(175,30){$\dots$} \put(275,30){\circle*{2}}
\put(227,35){$\mbox{}_{-2}$}\put(235,30){\circle*{3}}
\put(233,19){$\mbox{e}_{k+2}$}\put(195,30){\line(1,0){40}}
\put(263,35){$\mbox{}_{-1}$}\put(275,30){\circle*{3}}
\put(270,19){$\mbox{e}_{k+3}$}\put(235,30){\line(1,0){40}}
\put(275,70){\circle*{3}}
\put(280,67){$\mbox{b}_{1}$}\put(275,30){\line(0,1){40}}
\put(315,30){\circle*{3}}
\put(313,19){$\mbox{b}_{2}$}\put(275,30){\line(1,0){40}}
\put(200,-8){$\text{Fig.}\, 1$}
\end{picture} \vspace{0.5cm}

The graph $\Gamma(\widetilde B)$ of the curve germ $(B,o)$ of singularity type ${\bf A}_{2k}$, $k\geq 1$, is depicted on Fig. 2.

\begin{picture}(300,85)
\put(85,30){\circle*{3}}\put(80,35){$\mbox{}_{-2}$}
\put(82,19){$\mbox{e}_{2}$}\put(85,30){\line(1,0){40}}
\put(125,30){\circle*{3}}\put(120,35){$\mbox{}_{-2}$}
\put(122,19){$\mbox{e}_{3}$}\put(125,30){\line(1,0){40}}
\put(175,30){$\dots$} \put(275,30){\circle*{2}}
\put(227,35){$\mbox{}_{-3}$}\put(235,30){\circle*{3}}
\put(233,19){$\mbox{e}_{k+1}$}\put(195,30){\line(1,0){40}}
\put(263,35){$\mbox{}_{-1}$}\put(275,30){\circle*{3}}
\put(270,19){$\mbox{e}_{k+2}$}\put(235,30){\line(1,0){40}}
\put(275,70){\circle*{3}}
\put(280,67){$\mbox{b}_{1}$}\put(275,30){\line(0,1){40}}
\put(307,35){$\mbox{}_{-2}$}\put(315,30){\circle*{3}}
\put(313,19){$\mbox{e}_{k+3}$}\put(275,30){\line(1,0){40}}
\put(200,-8){$\text{Fig.}\, 2$}
\end{picture} \vspace{0.5cm}

The graph $\Gamma(\widetilde B)$ of the curve germ $(B,o)$ of singularity type ${\bf D}_{2k+2}$, $k\geq 1$, is depicted on Fig. 3.

\begin{picture}(300,85)
\put(85,30){\circle*{3}}
\put(82,19){$\mbox{b}_{1}$}\put(85,30){\line(1,0){40}}
\put(125,30){\circle*{3}}\put(120,35){$\mbox{}_{-2}$}
\put(122,19){$\mbox{e}_{4}$}\put(125,30){\line(1,0){40}}
\put(175,30){$\dots$} \put(275,30){\circle*{2}}
\put(227,35){$\mbox{}_{-2}$}\put(235,30){\circle*{3}}
\put(233,19){$\mbox{e}_{k+2}$}\put(195,30){\line(1,0){40}}
\put(263,35){$\mbox{}_{-1}$}\put(275,30){\circle*{3}}
\put(270,19){$\mbox{e}_{k+3}$}\put(235,30){\line(1,0){40}}
\put(275,70){\circle*{3}}
\put(280,67){$\mbox{b}_{2}$}\put(275,30){\line(0,1){40}}
\put(315,30){\circle*{3}}
\put(313,19){$\mbox{b}_{3}$}\put(275,30){\line(1,0){40}}
\put(200,-8){$\text{Fig.}\, 3$}
\end{picture} \vspace{0.5cm}

The graph $\Gamma(\widetilde B)$ of the curve germ $(B,o)$ of singularity type ${\bf D}_{2k+3}$, $k\geq 1$, is depicted on Fig. 4.

\begin{picture}(300,85)
\put(85,30){\circle*{3}}
\put(82,19){$\mbox{b}_{1}$}\put(85,30){\line(1,0){40}}
\put(125,30){\circle*{3}}\put(120,35){$\mbox{}_{-2}$}
\put(122,19){$\mbox{e}_{3}$}\put(125,30){\line(1,0){15}}
\put(150,30){$\dots$} \put(275,30){\circle*{2}}
\put(175,30){\line(1,0){15}}\put(193,30){\circle*{3}}\put(192,35){$\mbox{}_{-2}$}
\put(190,19){$\mbox{e}_{k+1}$}
\put(227,35){$\mbox{}_{-3}$}\put(235,30){\circle*{3}}
\put(233,19){$\mbox{e}_{k+2}$}\put(195,30){\line(1,0){40}}
\put(263,35){$\mbox{}_{-1}$}\put(275,30){\circle*{3}}
\put(270,19){$\mbox{e}_{k+3}$}\put(235,30){\line(1,0){40}}
\put(275,70){\circle*{3}}
\put(280,67){$\mbox{b}_{2}$}\put(275,30){\line(0,1){40}}
\put(307,35){$\mbox{}_{-2}$}\put(315,30){\circle*{3}}
\put(313,19){$\mbox{e}_{k+4}$}\put(275,30){\line(1,0){40}}
\put(200,-8){$\text{Fig.}\, 4$}
\end{picture} \vspace{0.5cm}

The graph $\Gamma(\widetilde B)$ of the curve germ $(B,o)$ of singularity type ${\bf E}_{6}$ is depicted on Fig. 5.

\begin{picture}(300,85)
\put(147,35){$\mbox{}_{-4}$}\put(155,30){\circle*{3}}
\put(153,19){$\mbox{e}_{2}$}
\put(183,35){$\mbox{}_{-1}$}\put(195,30){\circle*{3}}
\put(190,19){$\mbox{e}_{3}$}\put(155,30){\line(1,0){40}}
\put(195,70){\circle*{3}}
\put(200,67){$\mbox{b}_{1}$}\put(195,30){\line(0,1){40}}
\put(227,35){$\mbox{}_{-2}$}\put(235,30){\circle*{3}}
\put(233,19){$\mbox{e}_{4}$}\put(195,30){\line(1,0){40}}
\put(267,35){$\mbox{}_{-2}$}\put(275,30){\circle*{3}}
\put(273,19){$\mbox{e}_{5}$}\put(235,30){\line(1,0){40}}
\put(200,-8){$\text{Fig.}\, 5$}
\end{picture} \vspace{0.5cm}

The graph $\Gamma(\widetilde B)$ of the curve germ $(B,o)$ of singularity type ${\bf E}_{7}$ is depicted on Fig. 6.

\begin{picture}(300,85)
\put(127,35){$\mbox{}_{-3}$}\put(135,30){\circle*{3}}
\put(132,19){$\mbox{e}_{3}$}
\put(163,35){$\mbox{}_{-1}$}\put(175,30){\circle*{3}}
\put(170,19){$\mbox{e}_{4}$}\put(135,30){\line(1,0){40}}
\put(175,70){\circle*{3}}
\put(180,67){$\mbox{b}_{2}$}\put(175,30){\line(0,1){40}}
\put(207,35){$\mbox{}_{-2}$}\put(215,30){\circle*{3}}
\put(213,19){$\mbox{e}_{5}$}\put(175,30){\line(1,0){40}}
\put(255,30){\circle*{3}}
\put(253,19){$\mbox{b}_{1}$}\put(215,30){\line(1,0){40}}
\put(200,-8){$\text{Fig.}\, 6$}
\end{picture} \vspace{0.5cm}

The graph $\Gamma(\widetilde B)$ of the curve germ $(B,o)$ of singularity type ${\bf E}_{8}$ is depicted on Fig. 7.

\begin{picture}(300,85)
\put(127,35){$\mbox{}_{-3}$}\put(135,30){\circle*{3}}
\put(132,19){$\mbox{e}_{5}$}
\put(163,35){$\mbox{}_{-1}$}\put(175,30){\circle*{3}}
\put(170,19){$\mbox{e}_{4}$}\put(135,30){\line(1,0){40}}
\put(175,70){\circle*{3}}
\put(180,67){$\mbox{b}_{1}$}\put(175,30){\line(0,1){40}}
\put(207,35){$\mbox{}_{-2}$}\put(215,30){\circle*{3}}
\put(213,19){$\mbox{e}_{3}$}\put(175,30){\line(1,0){40}}
\put(247,35){$\mbox{}_{-3}$}\put(255,30){\circle*{3}}
\put(253,19){$\mbox{e}_{2}$}\put(215,30){\line(1,0){40}}
\put(200,-8){$\text{Fig.}\, 7$}
\end{picture} \vspace{0.5cm}

\begin{rem} \label{rem1} Note that in all graphs $\Gamma(\widetilde B)$ of a curve germs $(B,o)$ of $ADE$ singularity types {\rm (}except the singularity types ${\bf A}_0$ and ${\bf A}_1${\rm )} there is a single vertex $e$ of valency three {\rm (}denote the corresponding to it curve by $E${\rm )} and this vertex has weight $w=-1$. \end{rem}

\begin{prop} {\rm (\cite{K-rig}, Corollary 1)} \label{corol} Let $(B,o)$ be a curve germ having one of $ADE$ singularity types, $E\subset \sigma^{-1}(o)\subset \widetilde V$ the exceptional curve of the last blowup 
in the sequence of blowups $\sigma:\widetilde V\to V$ resolving the singular point of $(B,o)$, and $e$ an  element in $\pi_1^{loc}(B,o)$ represented by a simple loop around $E$. Then $e$ belongs to the center of $\pi_1^{loc}(B,o):= \pi_1(V\setminus B)\simeq  \pi_1(\widetilde V\setminus \widetilde B)$.
\end{prop}
\begin{prop} {\rm (\cite{K-rig}, Proposition 1)} \label{gener} Let $(B,o)$ be a curve germ having one of $ADE$ singularity types.  If the singularity type of $(B,o)$ is not $A_0$ or $A_1$, then $\pi_1^{loc}(B,o)$ is generated by  $e$ and the elements $\gamma_1,\gamma_2,\gamma_3$ corresponding to the vertices of $\Gamma(\widetilde B)$ connected by an edge with the vertex $e$ {\rm (if the singularity type of $(B,o)$ is  $A_1$, then $\pi_1^{loc}(B,o)$ is generated by  $b_1,b_2$ and $e$)}.
\end{prop}

Below, if  the singularity type of $(B,o)$ is ${\bf A}_{2n+1}$ or ${\bf D}_{2n+2}$, then we will identify the element $\gamma_1$ with $e_{n+2}$  (see Fig. 1 and Fig. 3); if the singularity type is ${\bf A}_{2n}$,  then we will identify 
$\gamma_1$ with $e_{n+1}$ and $\gamma_2$ with $e_{n+3}$ (see Fig. 2);  if the singularity type is ${\bf D}_{2n+3}$,  then we will identify 
$\gamma_1$ with $e_{n+2}$ and $\gamma_2$ with $e_{n+4}$ (see Fig. 4);    if the singularity type is ${\bf E}_{6}$,  then we will identify 
$\gamma_1$ with $e_{2}$ and $\gamma_2$ with $e_{4}$ (see Fig. 5);  if the singularity type is ${\bf E}_{7}$,  then we will identify 
$\gamma_1$ with $e_{5}$ and $\gamma_2$ with $e_{3}$ (see Fig. 6); and  if the singularity type is ${\bf E}_{8}$,  then we will identify 
$\gamma_1$ with $e_{5}$ and $\gamma_2$ with $e_{3}$ (see Fig. 7).

Let $\overline{\widetilde B\setminus E}$ be the closure of $\widetilde B\setminus E$ in $\widetilde V$.
\begin{df} \label{tail} If the singularity type of $(B,o)$ is not ${\bf A}_0$ or ${\bf A}_1$, then $\overline{\widetilde B\setminus E}$ is a disjoint union of three chains of curves which we call {\it trails} of $\widetilde B$. Denote by $\widetilde B_j$, $j=1,2,3$, the trail containing the curve for which $\gamma_j$ is represented by a loop around this curve. A trail is  {\it exceptional} {\rm (}resp., {\it completely exceptional}{\rm )} if it contains an exceptional curve of $\sigma$ {\rm (}resp., it contains only exceptional curves{\rm ).} \end{df}

Denote by $Z_e$ the subgroup of $\pi_1^{loc}(B,o)$ generated by $e$. The imbedding $i_1:\widetilde V\setminus \widetilde B\hookrightarrow \widetilde V\setminus (\widetilde B_1\cup\widetilde B_2\cup\widetilde B_3)$ induces an epimorphism $i_{1*}: \pi_1(\widetilde V\setminus \widetilde B)\to \pi_1(\widetilde V\setminus (\widetilde B_1\cup\widetilde B_2\cup\widetilde B_3))$ whose kernel is $Z_e$. It easily follows from Theorem \ref{Mum} that  $e=\gamma_1\gamma_2\gamma_3$ and
\begin{equation}\label{res1} \pi_1(\widetilde V\setminus (\widetilde B_1\cup\widetilde B_2\cup\widetilde B_3))\simeq \langle \widetilde{\gamma}_1\rangle * \langle \widetilde{\gamma}_2\rangle * \langle \widetilde{\gamma}_3\rangle / (\langle \widetilde{\gamma}_1\widetilde{\gamma}_2\widetilde{\gamma}_3\rangle)
\end{equation}
is the quotient group of free product of three cyclic groups $\langle \widetilde{\gamma}_j\rangle$, $j=1,2,3$, by the normal closure $(\langle \widetilde{\gamma}_1\widetilde{\gamma}_2\widetilde{\gamma}_3\rangle)$ of cyclic subgroup  generated by the product $\widetilde{\gamma}_1\widetilde{\gamma}_2\widetilde{\gamma}_3$, where
$\widetilde{\gamma}_j=i_{1*}(\gamma_j)$. It follows from Theorem \ref{Mum} that the group $\langle \widetilde{\gamma}_j\rangle$ is finite if and only if $\widetilde B_j$ is a union of exceptional curves $E_l$.

It is easy to see that the imbedding
$i_2: E\setminus (\widetilde B_1\cup\widetilde B_2\cup\widetilde B_3)\hookrightarrow \widetilde V\setminus  (\widetilde B_1\cup\widetilde B_2\cup\widetilde B_3)$
induces an epimorphism
\begin{equation}\label{res2} i_{2*}: \pi_1(E\setminus (\widetilde B_1\cup\widetilde B_2\cup\widetilde B_3),p)\to \pi_1(\widetilde V\setminus (\widetilde B_1\cup\widetilde B_2\cup\widetilde B_3),p)
\end{equation}
(here we assume that $p\in E\setminus (\widetilde B_1\cup\widetilde B_2\cup\widetilde B_3)$).
\begin{rem} \label{ident1} Note that if the singularity type of $(B,o)$ is $\bf D_4$ then $i_{2*}$ is an isomorphism. Therefore in this case we will identify the 
groups $\pi_1(E\setminus (\widetilde B_1\cup\widetilde B_2\cup\widetilde B_3),p)$ and
$\pi_1(\widetilde V\setminus (\widetilde B_1\cup\widetilde B_2\cup\widetilde B_3),p)$.\end{rem}

Denote by $P_j=E\ \cap \widetilde B_j$ and $\overline{\gamma}_j$ a loop around $P_j$ such that $i_{2*}(\overline{\gamma}_j)=\widetilde{\gamma}_j$.

\begin{df} \label{trailpi} If $\widetilde B_j$ is an exceptional trail, denote by $\widetilde B^0_j$ the union of the exceptional curves contained in $\widetilde B_j$ and by
$\widetilde{\pi}_j:=\pi_1(N_T\setminus \widetilde B_j)$, $\widetilde{\pi}_j^0:=\pi_1(N_T\setminus \widetilde B^0_j)$, where $N_T$ is a sufficiently small tubular neighbourhood of $\widetilde B_j$. \end{df}
\begin{rem} \label{tailcyc} It follows from Theorem {\rm \ref{Mum}} that $\widetilde{\pi}_j$ and $\widetilde{\pi}_j^0$ are cyclic groups. \end{rem}

\subsection{Cyclic quotients.}\label{H-Js} Let a cyclic group $\mu_m\simeq\mathbb Z_m$ of order $m$ act on a germ of a smooth surface $(U,o')$. Denote by $(W,o_1)=(U,o')/\mu_m$  the quotient space and  $\xi: (U,o')\to (W,o_1)$ the quotient map. By Cartan's Lemma, we can assume that $(U,o')$ is biholomorphic to the ball $\mathbb B_2=\{ (u_1,u_2)\in \mathbb C^2\mid |u_1|^2+|u_2|^2<1\}$ and the action of a generator $g$ of $\mu_m$ is given by
$$ g:(u_1,u_2)\mapsto (exp(2\pi  p_1i/m)u_1,exp(2\pi  p_2i/m)u_2),$$
where $p_j$, $j=1,2$, are some integers, $1\leq p_j\leq m$, $G.C.D.(m,p_1,p_2)=1$. Let $m=m_1m_2m_0$, $p_1=m_1t_1s$, $p_2=m_2t_2s$, where
$$G.C.D.(m_1t_1,m_2t_2)=G.C.D.(st_1,m_0)=G.C.D.(st_2,m_0)=1.$$
Then
$$ g^{m_1m_0}: (u_1,u_2)\mapsto  (exp(2\pi  p_1i/m_2)u_1,u_2), \,\,  g^{m_2m_0}: (u_1,u_2)\mapsto  (u_1,exp(2\pi  p_2i/m_1)u_2)$$
and the subgroup $\mu_{m_1m_2}\subset \mu_m$, generated by $g^{m_1m_0}$ and $g^{m_2m_0}$, is a cyclic group of order $m_1m_2$. The map $\xi$ can be decomposed into a composition of two maps, $\xi=\varphi\circ\vartheta_{m_1,m_2}$, where $\vartheta_{m_1,m_2}:(U,o')\to (X,\widetilde o)$ is the quotient map defined by the action of $\mu_{m_1m_2}$ on $(U,o')$ and $\varphi:(W,o_1)\to (V,o)$  is the quotient map defined by the action of the quotient group $\mu_m/\mu_{k_1k_2}\simeq \mu_n$ of order $n$ on $(W,o_1)$.

It is easy to see that $(X,\widetilde o)$ is a germ of a smooth surface,
$$(X,\widetilde o)\simeq \mathbb B_2=\{ (x_1,x_2)\in \mathbb C^2\mid |x_1|^2+|x_2|^2<1\},$$ and  the map $\vartheta_{m_1,m_2}$ is given by $x_1=u_1^{m_2}$, $x_2=u_2^{m_1}$. The image $\overline g$ in $\mu_{m_0}$ of the generator $g\in\mu_m$ acts on $(X,\widetilde o)$ as follows:
$$\overline g:(x_1,x_2)\mapsto (exp(2\pi  st_1i/n)w_1,exp(2\pi  st_2i/n)w_2).$$ There is an integer $r$ such that
$rst_2\equiv 1\,\, \text{mod}\, m_0$ and $rst_1\equiv q\,\, \text{mod}\, m_0$,  where $1\leq m_0$, since $G.C.D.(st_j,m_0)=1$ for $j=1,2$. Therefore
$$\overline g^r:(x_1,x_2)\mapsto (exp(2\pi qi/m_0)x_1,exp(2\pi i/m_0)x_2)$$ and it is easy to show that $(W,o_1)$
 is the normalization of the germ of singularity in $\mathbb B_3=\{ (z_1,z_2,z_3)\in \mathbb C^3\mid |z_1|^2+|z_2|^2+|z_3|^2<1\}$ given by $z_3^n=z_1z_2^{n-q}$, where $z_3=x_1x_2^{n-q}$, $z_1=x_1^{m_0}$, and $x_2=x_2^{m_0}$, that is, the germ $(W,o_1)$ has so called {\it Hirzebruch-Jung singularity  type} $A_{m_0,q}$.

 The map $\vartheta_{m_1,m_2}$ is branched in $L_1=\{ x_1=0\}$ (if $m_2>0$) and  $L_2=\{ x_2=0\}$  (if $m_1>0$), and $\varphi$ is unramified outside $\widetilde o$ (in the future we will denote the map $\varphi$ by $\theta_{m_0,q}$). Therefore $\pi_1(W\setminus o_1)\simeq \mu_{m_0}$ and $\theta_{m_0,q}: X\setminus \widetilde o\to W\setminus o_1$ is the universal unramified cover.

Let $\tau: \widetilde W\to (W,o_1)$ be the minimal resolution of the singular point $o_1\subset W$. Denote by $B_j=\tau^{-1}(\theta_{m_0,q}(L_j))$, $j=1,2$, the proper inverse image of $\theta_{m_0,q}(L_j)$ and let $\tau^{-1}(o_1)=\bigcup_{j=1}^kE_j$. It is well-known (see, for example, \cite{B}, Capter III.5) that $E_j$ are rational curves and up to renumbering of $E_j$, the dual weighted graph $\Gamma(\widetilde B)$ of the curve $\widetilde B=(B_1\cup B_2)\cup(\bigcup_{j=1}^kE_j)$ is a chain, that is it has the following form:

\begin{picture}(300,60)
\put(99,38){$B_1$}\put(105,30){\circle{5}}
\put(145,30){\circle*{5}}\put(141,38){$E_1$}
\put(141,20){$\omega_1$}\put(108,30){\line(1,0){40}}
\put(145,30){\line(1,0){40}}\put(215,30){\line(1,0){40}} \put(255,30){\circle*{5}}
\put(195,30){$.\, .\, .$}
\put(253,20){$\omega_k$} \put(251,38){$E_k$} \put(257,30){\line(1,0){40}} \put(299,30){\circle{5}}
\put(295,38){$B_2$}
\put(190,-8){$\text{Fig.}\, 8$}
\end{picture} \vspace{0.5cm}

where $\omega_j=-(E_j^2)_{\widetilde W}$ satisfy the following equality:
\begin{equation}\label{H-J}
\frac{m_0}{q}=\omega_1-\frac{1}{\omega_2-\frac{1}{\omega_3- \frac{1}{\dots -\frac{1}{\omega_k}}}}:\stackrel{def}{=}[w_1;w_2,\dots,w_k].\end{equation}

Conversely, if $\tau:\widetilde W\to (W,o_1)$ is the minimal resolution of normal singularity such that $\tau^{-1}(o_1)=\bigcup_{j=1}^kE_j$ is a chain of rational curves 
(see Fig. 8), then $(W,o_1)$ is a Hirzebruch-Jung singularity of type $A_{m_0,q}$, where $m_0$ and $q$ can be found using equality (\ref{H-J}).

\begin{rem} \label{unique} A representation of the singularity $(W,o_1)$ of type $A_{m_0,q}$ as a cyclic quotient singularity is uniquely  defined  by a choice of divisors $B_1$ and $B_2$
{\rm (}see Fig. 8{\rm )} in $\widetilde W$ {\rm (}\cite{B}, Chapter III.5{\rm )}.
\end{rem}

Note that if we renumber the curves $E_j$ and the weights $\omega_j$ as follows: $E'_j:=E_{k-j+1}$ and $\omega'_j:=\omega_{k-j+1}$, and substitute the new $\omega'_j$ instead of the old $\omega_j$ in the right side of (\ref{H-J}) then we obtain $\frac{m_0}{q'}$ in the left side of (\ref{H-J}) with $q'$ such that $qq'\equiv 1\, {\text{mod}}\, m_0$ (see \cite{B}, Chapter III.5). In particular, the singularities of types $A_{m_0,q}$ and $A_{m_0,q'}$ are the same singularity.

\begin{rem}\label{gener} In notations used in Theorem {\rm \ref{Mum}}, it easily follows from Theorem {\rm \ref{Mum}} that the group $\pi_1(\widetilde W\setminus (B_1\cup B_2\cup(\bigcup_{j=1}^kE_j)))$ is generated by the elements $b_1$ and $e_1$, and $\pi_1(\widetilde W\setminus ( B_2\cup(\bigcup_{j=1}^kE_j)))$ is the free group $\mathbb F_1$ generated by $e_1$.
\end{rem}

For a singularity $(W,o_1)$ of type $A_{m_0,q}$ we have $$\pi_1(W\setminus \{ o_1\})=\pi_1(\widetilde W\setminus (\bigcup_{j=1}^k E_j))\simeq \mu_{m_0},$$  since $\theta_{m_0,q}: X\setminus \{\widetilde o\} \to W\setminus \{o_1\}$ is the universal cover.

\begin{lem} \label{A_n0} If $[\omega_1,\omega_2,\dots,\omega_{k}]=[2,\dots,2]$, $k\geq 1$, then $(W,o_1)$ has the singularity type  $A_{k+1,k}$ and, in particular, $\pi_1(W\setminus \{o_1\})\simeq \mu_{k+1}$.
\end{lem}
\proof
We have
\begin{equation} \label{A_n}[\underbrace{2;2,\dots,2}_{k}]=\frac{k+1}{k}. \end{equation}
Note that the singularity types ${\bf A}_{k}$ and $A_{k+1,k}$ are the same singularity type. \qed

\begin{lem} \label{3A_n} If $[\omega_1,\omega_2,\dots,\omega_{k+1}]=[n,2,\dots,2]$, $k\geq 0$, then $(W,o_1)$ has the singularity type  $A_{n(k+1)-k,k+1}$  and, in particular, $\pi_1(W\setminus \{o_1\})\simeq \mu_{n(k+1)-k}$.\end{lem}
\proof We have
\begin{equation}\label{rel} [\omega_1;\omega_2,\dots,\omega_{k+1}]=\omega_1-\frac{1}{[\omega_2;\dots,\omega_{k+1}]}.\end{equation}
Therefore $[n;2,\dots,2]=n-\frac{k}{k+1}=\frac{n(k+1)-k}{k+1}$. \qed \\

Denote by $\mathbb D^2_{r_1,r_2}=\{ (y_1,y_2)\in \mathbb C^2\mid |y_1|<r_1, |y_2|< r_2 \}$ a bidisk in $\mathbb C^2$, where
$(r_1,r_2)\in \mathbb R_+^2$.
Let $L_{y_j}=\{ y_j=0\}\subset \mathbb D^2_{r_1,r_2}$, $j=1,2$,  be the coordinates axis in $\mathbb D^2_{r_1,r_2}$.

The following Lemma is a direct consequence of Theorem 5.1 in \cite{B}, Chapter III.
\begin{lem}\label{cyc2} Let $Z$ be an irreducible germ of normal surface and $\xi:Z\to\mathbb D^2_{(r_1,r_2)}(y_1,y_2)$ a cyclic cover of degree $n$ branched in $L_{y_1}\cup L_{y_2}$. Then $n=n_1n_2n_3$ for some $n_1\geq 1$, $n_2\geq 1$, $n_3\geq 1$ such that $G.C.D.(n_1,n_2)=1$ and $\xi$ is ramified over $L_j$ with multiplicity $n_jn_3$, $j=1,2$,  and if $n_3>1$ then the singularity type of $Z$ over the origin $(0,0)\in \mathbb D^2_{r_1,r_2}$ is $A_{n_3,q}$ for some $q$, $G.C.D.(n_3,q)=1$, and 
$Z$ is a germ of smooth surface if $n_3=1$.
\end{lem}

Let a germ  $(W,o_1)$ of normal surface have the singularity type $A_{n,q}$, $\tau: \widetilde W\to (W,o_1)$ the minimal resolution of the singular point $o_1\in W$, $\tau^{-1}(o_1)=\bigcup_{j=1}^kE_j$ the chain of rational curves, $(E_j)^2=-\omega_j$, and let $\widetilde B\subset \widetilde W$ be a curve whose dual graph $\Gamma(\widetilde B)$ is depicted on Fig. 8.

Let $m$ be a divisor of $n$, $n=mk$. Denote by $\varphi_m :(X_m,\widetilde o)\to (W,o_1)$ the cyclic cover of degree $m$ defined by the natural epimorphism
$\varphi_{m*}: \pi_1(W\setminus \{ o_1\})\simeq \mu_n\twoheadrightarrow \mu_m$. The cover $\varphi_m$ is unramified outside $\widetilde o$ and $(X_m,\widetilde o)$ is a normal variety having the singularity of type $A_{k,q'}$ for some $q'$ if $k>1$, and $(X_m,\widetilde o)$ is a germ of smooth surface if $k=1$.

Consider the following commutative diagram
\begin{equation}\label{diag}
\begin{CD} \overline X_m @>\varrho>> \widetilde X_m @>\rho>> (X_m,\widetilde o) \\
@. @V\widetilde{\varphi}_m VV @VV\varphi_m V \\
\phantom{\overline W} @. \widetilde W @>\tau>> (W,o_1) \end{CD}
\end{equation}
in which $\widetilde X_m$ is the normalization of fibre product $\widetilde W\times_{(W,o_1)} (X_m,\widetilde o)$ and $\varrho:\overline X_m\to \widetilde X_m$ is the minimal resolution of singular points of $\widetilde X_m$. It follows from Lemma \ref{cyc2} that $\widetilde X_m$ can have singular points (and their singularity types are $A_{m',q'}$ with some divisors $m'$ of $m$) only over the intersection points of neighboring exceptional curves $E_j$ and $E_{j+1}$ of $\tau$.

Denote by $\psi:=\rho\circ\varrho$ the composition of maps $\rho$ and $\varrho$. Note that $\psi$ is a resolution of singular point $\widetilde o\in X_m$. The map $\psi$ is decomposed into a composition of two maps, $\psi=\varsigma\circ\sigma$, where $\varsigma:\overline X_{m,min}\to (X_m,\widetilde o)$ is the minimal resolution of the singular point $\widetilde o$ if $m<n$ and  $\sigma:\overline X_m\to \overline X_{m,min}$ is a composition of $\sigma$-processes, $\sigma=\sigma_l\circ \dots\circ \sigma_1$ ($\psi=\sigma$ if $m=n$).

Let $\widetilde B=(B_1\cup B_2\cup(\bigcup_{j=1}^kE_j\subset \widetilde W$ be the union of curves and curve germs whose dual weighted graph is depicted in Fig. 8. Denote by the same letter $C_j$  the proper inverse image $(\tau\circ\widetilde{\varphi}_m)^{-1}(B_j)$ of the germ $B_j$, $j=1,2$, and the proper inverse images  $(\sigma_l\circ\dots\circ\sigma_{l-s})^{-1}(C_j)$ for $1\leq s\leq l$. Denote by $\Delta_m(n,q)$ the number of $\sigma$-processes from the set $\{\sigma_1,\dots\sigma_{l}\}$ which blowup a point lying in $C_1$ and call it the {\it $m$-th supplement for the singularity type} $A_{n,q}$.

\begin{lem}\label{delta1} Let a germ $(W,o_1)$ have the singularity type $A_{n,n-1}$, $n=mk$, and  $\varphi_m :(Z_m,\widetilde o)\to (W, o_1)$ the cyclic cover of degree $m$ defined by the natural epimorphism
$\varphi_{m*}: \pi_1(W\setminus \{ o_1\})\simeq \mu_n\twoheadrightarrow \mu_m$. Then
\begin{itemize}
\item[$(i)$] the singularity type of $(Z_m,\widetilde o)$ is $A_{k,k-1}$\,\footnote{By definition, $(Z_m,\widetilde o)$ has the singularity type $A_{1,0}$ (when $k=1$) means that the point $\widetilde o$ is a smooth point of $Z_m$.},
\item[$(ii)$] $\Delta_m(n,n-1)=m-1$.
\end{itemize}
\end{lem}

\proof To prove $(ii)$,  let us consider the quadric  $Q=\mathbb P^1\times \mathbb P^1$ and let $S_1$, $S_2$ be two fibres of the projection $\text{pr}_2:Q\to \mathbb P^1$ to the second factor and $L$ a fibre of the projection $\text{pr}_1:Q\to \mathbb P^1$ to the first factor. Consider the following diagram
\begin{equation}\label{diag1}
\begin{CD} \overline X_m @>\varrho>> \widetilde X_m @>\rho>> X_m \\
@. @V\widetilde{\varphi}_m VV @VV\varphi_m V \\
\phantom{\overline V} @. \widetilde Q @>\tau>> Q \end{CD}
\end{equation}
in which
\begin{itemize}
\item[$1)$]  $\varphi_m$ is defined by epimorphism $\varphi_{m*}:\pi_1(Q\setminus (S_1\cup S_2))\simeq \mathbb Z\to\mu_m\subset\mathbb S_m$,
\item[$2)$] $\tau=\tau_1\circ\dots\circ\tau_n$ is the composition of $n$ blowups of the point $p=S_1\cap L\in S_1$,
\item[$3)$] $\widetilde X_m$ is the normalization of the fibre product $\widetilde Q\times_Q X_m$,
\item[$4)$] $\varrho:\overline X_m\to \widetilde X_m$ is the minimal resolution of the singular points of $\widetilde X_m$.
\end{itemize}

Denote by $E_j\subset \widetilde Q$, $j=1,\dots,n-1$, the proper inverse image of the  exceptional curve of blowup $\tau_j$, $B_1\subset \widetilde Q$ the proper inverse image of the  the fibre $L$, and $B_2\subset \widetilde Q$ the  exceptional curve of blowup $\tau_n$. We have
$$ (B_1^2)_{\widetilde Q}=(B_2^2)_{\widetilde Q}=-1,\quad (E_j^2)_{\widetilde Q}=-2\,\,\, \text{for}\,\, j=1,\dots,n-1$$
and the dual graph $\Gamma(\widetilde B)$ of $\widetilde B=B_1\cup B_2\cup(\bigcup_{j=1}^{n-1}E_j)$ is shown in Fig. 8  (in which $\omega_j=-2$ for $j=1,\dots, n-1$). Therefore we can identify $\widetilde W$ with a tubular neighbourhood of $\bigcup_{j=1}^{n-1}E_j$ and $\widetilde Z_m$ with $\widetilde{\varphi}_m(\widetilde W)$. Note that the dual weighted graph $\Gamma(\bigcup_{j=1}^{n-1}E_j)$ of the curve $\bigcup_{j=1}^{n-1}E_j$ is a {\it central-symmetric graph}, that is, the weights $\omega_j$ of the vertices  $e_j$ satisfy the following relation: $\omega_j=\omega_{n-j}$.

Obviously, $X_m\simeq \mathbb P^1\times\mathbb P^1$, where $\varphi_m^{-1}(S_1)$, $\varphi_m^{-1}(S_2)$ are two fibres of the projection $\text{pr}_2:X_m\to\mathbb P^1$ to the second factor and $\varphi_m^{-1}(L):=F$ is a fibre of the projection $\text{pr}_1:X_m\to \mathbb P^1$ to the first factor.

The fundamental group $\pi_1(Q\setminus (S_1\cup S_2))\simeq \pi_1(\widetilde Q\setminus (\tau^{-1}(S_1\cup S_2)\cup B_2\cup (\bigcup_{j=1}^{n-1}E_j))$ is generated by the element $\gamma$ represented by a simple loop around the curve $S_1$. Denote by $e_j$, $j=1,\dots,n-1$, and $b_2$ the elements of $\pi_1(\widetilde Q\setminus (\tau^{-1}(S_1\cup S_2)\cup B_2\cup (\bigcup_{j=1}^{n-1}E_j))$ represented by simple loops respectively around $E_j$ and $B_2$. It follows from Lemma \ref{sigma}  that
\begin{equation}\label{chain} e_j=\gamma^j\,\,\, \text{for}\,\,\, j=1,\dots, n-1\,\,\, \text{and}\,\,\, b_2=\gamma^n.\end{equation}

In the beginning, consider the case $m=n$. The element $g_1=\varphi_{n*}(\gamma)$ is a generator of $\mu_n\subset \mathbb S_n$ and
\begin{equation}\label{symmet1} \varphi_{n*}(e_j)=g_1^{j}, \end{equation}
in particular,  $\varphi_{n*}(e_{n-1})=g_1^{n-1}:=g_2$ is a generator of $\mu_n$ and $\varphi_{n*}(b_2)=g_1^n=\text{id}$. Therefore $\varphi_n$ is not branched in $B_2$ and $\widetilde X_n$ is smooth in a neighbourhood of $\varphi_n^{-1}(B_1)$ and $\varphi_n^{-1}(B_2)$. The restriction of $\widetilde{\varphi}_n$ to $\widetilde Z_n\subset \widetilde X_n$ is defined by monodromy homomorphism $\varphi_{n*}:\pi_1(\widetilde W\setminus \bigcup_{j=1}^{n-1}E_j)\to\mu_n\subset\mathbb S_n$ sending $e_j$ to $g_1^j$. It easily follows from Theorem \ref{Mum} that $\pi_1(\widetilde W\setminus \bigcup_{j=1}^{n-1}E_j)$ is generated by $e_{n-1}$ and $e_j=e_{n-1}^{n-j}$. Therefore if we denote $e'_{j}=e_{n-j}$ then \begin{equation}\label{symmet2} \varphi_{n*}(e'_j)=g_2^{j}. \end{equation}
The inverse image $(\varphi_n\circ\varrho)^{-1}(\bigcup_{j=1}^{n-1}E_j)=\bigcup_{j=1}^{N}\overline E_j\subset \overline X_n$ is the chain of rational curves which can be contracted to a smooth point, the curves $\overline B_j=(\varphi_n\circ\varrho)^{-1}(B_j)$, $j=1,2$, are rational curves,
$$(\overline B_j^2)_{\overline X_n}=\deg \widetilde{\varphi}_n\cdot (B_j^2)_{\widetilde X_n}=-n,\,\,\, (\overline B_1,\overline E_1)_{\overline X_n}=(\overline B_2,\overline E_N)_{\overline X_n}=1,$$
and $(\varphi_n\circ\rho\circ\varrho)^{-1}(L)=\overline B_1\cup\overline B_2\cup(\bigcup_{j=1}^{N}\overline E_j)$ is a fibre of the ruled surface $\overline X_n$.

It follows from central symmetry of the graph  $\Gamma(\bigcup_{j=1}^{n-1}E_j)$  and (\ref{symmet1}), (\ref{symmet2}) that the dual weighted graph $\Gamma(\bigcup_{j=1}^{N}\overline E_j)$ is also central-symmetric. Therefore if $\overline E_j$ is a $(-1)$-curve, then  $\overline E_{N-j+1}$ is also a $(-1)$-curve and the curves $\overline E_j$, $\overline E_{N-j+1}$ can be contracted to points simultaneously. After consistently contractions of all such pairs of $(-1)$-curves and the contraction of the curve $\overline E_{K+1}$ in the last step (it is easy to see that $N=2K+1$ must be an odd number and the central curve $\overline E_{K+1}$ is contracted in the last step of contractions) we obtain that the images of $\overline B_1$ and $\overline B_2$ are $(-1)$-curves, since the union of these images is a fibre of ruled structure. Therefore $\Delta_{n}(n,n-1)=n-1$, since $(B_1^2)_{\overline X_n}=(B_2^2)_{\overline X_n}=-n$.

Consider the case when $m<n$. The element $g_1=\varphi_{m*}(\gamma)$ is a generator of $\mu_m\subset\mathbb S_m$. It follows from (\ref{chain}) that
$ \varphi_{m*}(e_{j+lm})=g_1^j$ for $ j=1,\dots,m-1$, $l=0,\dots,k-1$ and $\varphi_{m*}(e_{lm})=\varphi_{m*}(b_2)=\text{id}$ for $l=1,\dots,k-1$.
Therefore $\widetilde{\varphi}_m$ is not branched in $E_{lm}$, $l=1,\dots, k-1$, and in $B_1$, $B_2$. Hence, $\widetilde X_m$ is smooth in a neighbourhood of
$\widetilde{\varphi}_m^{-1}(B_1\cup B_2\cup(\bigcup_{l=1}^{k-1}E_{lm}))$ and
$$(\widetilde{\varphi}_m^{-1}(E_{lm}),\widetilde{\varphi}_m^{-1}(E_{lm}))_{\widetilde X_m}=-2m,\qquad (\widetilde{\varphi}_m^{-1}(B_j),\widetilde{\varphi}_m^{-1}(B_j))_{\widetilde X_m}=-m$$for $l=1,\dots,k-1$ and $j=1,2$. For each $l=0,\dots,k-1$, the union $\bigcup_{j=1}^{m-1}\widetilde{\varphi}_m^{-1}(E_{j+lm})$ can be contracted to a smooth point and similar to the case $m=n$, it is easy to see, first, $\Delta_m(n,n-1)=\Delta_m(m,m-1)=m-1$, second, after contraction of curves $\bigcup_{l=0}^{k-1}\bigcup_{j=1}^{m-1}\widetilde{\varphi}_m^{-1}(E_{j+lm})$ the images of curves $\widetilde{\varphi}_m^{-1}(E_{lm})$ form a chain consisting of $k-1$ $(-2)$-curves, that is the singularity type of $(Z_m,\widetilde o)$ is $A_{k,k-1}$. \qed \\

\section{Description of $\beta^{-1}(f)$: general case} \label{section2}
\subsection{Necessary conditions.} \label{neces} In this Section, we use notations of the previous Section.

Consider a rigid germ of cover $F:(U,o')\to (V,o)$ branched in a germ $(B,o)$ having one of $ADE$  singularity types at the point $o$, $d=\deg F$, and let $F_*:\pi_1^{loc}(B,o)=\pi_1(V\setminus B, p)\twoheadrightarrow G_F\subset \mathbb S_d$ be its monodromy homomorphism. Remind that the symmetric group $\mathbb S_d$ acts (from the right) on the fibre $F^{-1}(p)=\{ q_1,\dots, q_d\}$ and the monodromy group $G_F$ is a transitive subgroup of $\mathbb S_d$. Denote by $G_F^1$ the subgroup of $G_F$ leaving fixed the point $q_1$. Then the action of $G_F$ on $F^{-1}(p)$ can be identified with the action of $G_F$ on the set of right cosets of the subgroup $G^1_F$.

By Proposition \ref{corol}, the cyclic group $F_*(Z_e) \subset G_F$, generated by $F_*(e)$, is a central subgroup of $G_F$ and by Proposition 13 in \cite{K-rig}, the group $F_*(Z_e)$ acts on $(U,o')$. Denote by $H_1:(U,o')\to (W,o_1)=(U,o')/F_*(Z_e)$ the quotient map, $\deg H_1=m=|F_*(Z_e)|$, where $(W,o_1)$ is a germ of normal surface.
By Proposition 13 in \cite{K-rig}, there is a finite map $H_2:(W,o_1)\to (V,o)$ such that
$F=H_2\circ H_1$, $\deg H_2=n=\frac{d}{m}$. The monodromy group $G_{H_1}\subset\mathbb S_m$ of $H_1$ is isomorphic to $F_*(Z_e)$ and by Remark 2 in \cite{K-rig}, the monodromy group $G_{H_2}\subset\mathbb S_n$ of $H_2$ is isomorphic to $G_F/N$, where $N$ is the maximal normal subgroup of $G_F$ contained in $G^1_F\times F_*(Z_e)\subset G_F$.

Denote by $\widetilde W$ the normalization of fibre product $\widetilde V\times_{(V,o)} (W,o_1)$, where $\sigma:\widetilde V\to (V,o)$  is the minimal resolution of singular point $o\in V$ of the curve germ $(B,o)$ (if $(B,o)$ has the singularity of type ${\bf A}_0$ or ${\bf A}_1$ then $\sigma$ consists of the single blowup)
and let $\widetilde H_2:\widetilde W\to \widetilde V$ and $\varsigma:\widetilde W\to (W,o_1)$ be the projections to the factors. In addition, denote by $\widetilde U$ the normalization of fibre product $\widetilde W\times_{(W,o_1)} (U,o')$ and let $\widetilde H_1:\widetilde U\to \widetilde W$ and $\tau:\widetilde U\to (U,o')$ be the projections to the factors. The group $G_{H_1}$ acts on $\widetilde U$ and $\widetilde H_1$ is also the quotient map.

Denote by $C=\widetilde H_2^{-1}(E)$  the proper inverse image of the exceptional curve $E$ of the last $\sigma$-process and let $f:C\to E$, $\deg f=\deg \widetilde H_2=n$, be the restriction of $\widetilde H_2$ to $C$ (by definition, $f=\beta(F)$). The group $F_*(Z_e)$ acts on $\widetilde U$ and it is easy to see that $\widetilde H_1^{-1}(C)$ is an irreducible curve. Therefore the curve $C$ is contained in the branch locus of $\widetilde H_1$ and $\widetilde H_1$ is branched in $C$ with multiplicity $m=\deg \widetilde H_1$.  For the same reason, the branch locus of $\widetilde H_2$ is contained in $\widetilde B\setminus E$, where $\widetilde B=\sigma^{-1}(B)$ is the inverse image of the germ $(B,o)$.
The dual graph of $\widetilde B$ is depicted in one of Fig. 1 -- 7.

Let $\widetilde B^0=\sigma^{-1}(o)$ and let $\widetilde B_j\subset \widetilde B$, $j=1,2,3$ be the trails of $\widetilde B$ (see Definition \ref{tail}). It follows from Remark \ref{tailcyc} and Lemma \ref{cyc2} that $\widetilde W$ can have singular points (and their singularity types are $A_{k',q'}$ with some divisors $k'$ of $n$) only over the intersection points of neighboring irreducible components of the trails of $\widetilde B$, since $\widetilde H_2$ is not branched in $E$.
Denote by $\varsigma_r:\overline W\to \widetilde W$ a resolution of the singular points of $\widetilde W$. Then $\varsigma\circ\varsigma_r:\overline W\to (W,o_1)$ is a resolution of the  singular point $o_1$ of $(W,o_1)$.

Denote by $\varsigma_1:\overline W\to \overline W_m$  a holomorphic bimeromorphic map, where $\overline W_m$ is a smooth surface and $\varsigma_1$ contacts to points the maximal number of irreducible components belonging to $(\widetilde H_1\circ \varsigma_r)^{-1}(\bigcup_{j=1}^3(\widetilde B_j)$. Then $\varsigma\circ\varsigma_r= \varsigma_m\circ \varsigma_1$, where $\varsigma_m:\overline W_{m}\to (W,o_1)$  is also a resolution of the singular point $o_1$ of $(W,o_1)$.

 Put $\overline B^0=\varsigma_m^{-1}(o_1)$. The composition $\overline H_1=\varsigma_m^{-1}\circ H_1:(U,o')\to \overline W_{m}$ is a meromorphic map such that the finite cover $\overline H_1:U\setminus \{ o'\} \to \overline W_{m}\setminus \overline B^0$ is naturally isomorphic to the cover $H_1:U\setminus \{ o'\}\to W\setminus \{o_1\}$. Note that  $\overline C=\varsigma_1\circ \varsigma_r^{-1}(C)$ is a component of $\overline B^0$.

The cover $H_1$ is branched at most in two irreducible curve germs in $(W,o_1)$ (see Subsection \ref{H-Js}). Denote by $\overline B=\overline B_1\cup\overline B_2\cup\overline B^0\subset \overline W_{m}$ the inverse image of these germs (of course, one of $\overline B_j$ or both of them can be empty). The dual graph of $\overline B$ is a chain similar to one depicted in Fig. 8.

As a result, we have the following commutative diagram.

\begin{picture}(300,150)
\put(135,130){$\widetilde U$}
\put(155,135){\vector(1,0){120}}
\put(215,140){$\tau$}
\put(280,130){$U\ni o'$}

\put(278,128){\vector(-4,-3){28}}
\put(245,118){$\overline H_1$}

\put(280,60){$W\ni o_1$}
\put(135,60){$\widetilde W$}
\put(155,65){\vector(1,0){120}}
\put(215,55){$\varsigma$}
\put(215,105){$\varsigma_1$}
\put(171,78){$\varsigma_{r}$}
\put(244,78){$\varsigma_{m}$}
\put(385,60){$(*)$}

\put(140,125){\vector(0,-1){50}}
\put(285,125){\vector(0,-1){50}}
\put(125,100){$\widetilde H_1$}
\put(289,100){$H_1$}

\put(180,95){$\overline W$}
\put(195,100){\vector(1,0){35}}
\put(235,95){$\overline W_{m}$}
\put(180,93){\vector(-4,-3){30}}
\put(245,93){\vector(4,-3){30}}

\put(75,60){$\phantom{aaa}$} \put(60,60){$\mathbb P^1\simeq C$}
\put(113,60){$\subset$}

\put(60,0){$\mathbb P^1\simeq E$}
\put(118,0){$\subset$}
\put(135,0){$\widetilde V$}
\put(157,5){\vector(1,0){120}}
\put(215,-5){$\sigma$}
\put(280,0){$V\ni o$}
\put(78,35){$f$}
\put(91,54){\vector(0,-1){40}}
\put(143,35){$\widetilde H_2$}
\put(140,54){\vector(0,-1){40}}
\put(288,35){$H_2$}
\put(285,54){\vector(0,-1){40}}
\end{picture}\vspace{0.5cm}

The maps $\widetilde H_1$ and $H_1$ in diagram $(*)$ are quotient maps under the action of cyclic group. Therefore  $H_1$ is a composition of two maps, $H_1=\theta_{n',q}\circ \vartheta_{m_1,m_2}$ (see Subsection \ref{H-Js}), where $m=n'm_1m_2$ and $G.C.D.(m_1,m_2)=1$, and hence
we obtain the following commutative diagram
\begin{equation}\label{diag-D_4} \begin{CD} \widetilde H_2\circ \widetilde H_1:\phantom{a} @. \widetilde{U} @>{\widetilde{\vartheta}_{m_1,m_2}}>> \widetilde{X} @>{\widetilde{\theta}_{n',q}}>>
 \widetilde W @>{\widetilde H_2}>>  \widetilde V \\
@. @ V\tau VV @ V\overline{\varsigma} VV @VV\varsigma V @VV\sigma V \\
H_2\circ H_1:\phantom{a} @. (U,o') @>>{\vartheta}_{m_1,m_2}>  (X,\widetilde o) @>>{\theta}_{n',q}>  (W,o_1) @ >>H_2> (V,o). \end{CD}\end{equation}

The monodromy homomorphism $f_*$ is the composition of two homomorphisms, $f_*=H_{2*}\circ i_{2*}$ (see (\ref{res1}) and (\ref{res2})).
\begin{rem}\label{ident2} In view of Remark {\rm \ref{ident1}} we will identify the monodromy homomorphism $\widetilde H_{2*}$ with monodromy homomorphism $\beta(F)_*=f_*$  in the case when $(B,o)$ has the singularity type $\bf D_4$.
\end{rem}

The monodromy group $G_f$ of the Belyi pair $f$ is isomorphic to $G_{H_2}=G_F/N\subset \mathbb S_n$ generated by $f_*(\overline{\gamma}_j)\in \mathbb S_n$, $j=1,2,3$, where $\overline{\gamma}_j\in \pi_1(E\setminus (\widetilde B_1\cup\widetilde B_2\cup\widetilde B_3),p)$ are elements defined in \S \ref{graphade}. Let
$$T_c(\widetilde H_2)=T_c(f)=\{c_1,c_2,c_3\},\,\,\,\, c_j=(n_{j,1},\dots, n_{j,k_j}),\,\, n=\sum_{l=1}^{k_j}n_{j,l}$$
be the set of cycle types of $\widetilde H_{2*}(\widetilde {\gamma}_j)=f_*(\overline{\gamma}_j)$.
For each $j=1,2,3$ the inverse image $\widetilde H_2^{-1}(\widetilde B_j)$ of a tail $\widetilde B_j$ is the disjoint union of $k_j$ connected components,
$\widetilde H_2^{-1}(\widetilde B_j)=\bigsqcup_{l=1}^{k_j}\widetilde B_{j,l}$. Properties of cyclic covers, described in Subsection \ref{H-Js}, imply the following {\it contractibility  condition}:\newline
$\phantom{aa}$ {\it $\widetilde B_{j,l}\cap \widetilde H_2^{-1}(\widetilde B_j^{0})$ in $\widetilde W$ can be contracted to a nonsingular point if and only if the \newline  $\phantom{aa}$ order $|\widetilde{\pi}^0_j|$\, of the\, group\, $\widetilde{\pi}^0_j$ is a divisor of\, $n_{j,l}$ if $\widetilde B_j$\, is an exceptional trail, and \newline $\phantom{aa}$ $|\widetilde{\pi}^0_j|=n_{j,l}$ if $\widetilde B_j$ is a completely exceptional trail, } \newline since $\widetilde B_j^0$ can be contracted to a point of singularity type $A_{|\widetilde{\pi}^0_j|,q}$.
Note that, by the same reason, $n_{j,l}$ are divisors of $|\widetilde{\pi}^0_j|$ if $\widetilde B_j$ is a completely exceptional trail.

Denote by $r_j(H_{2*})$ the number of cycles in the permutation $f_*(\overline{\gamma}_j)$ whose lengths do not satisfy the contractibility condition if $\widetilde B_j$ is exceptional trail and put $r_j(H_{2*})=0$ if $\widetilde B_j$ is not an exceptional trail.
We obtain that
 \begin{equation}\label{nes} r_1(\widetilde H_{2*})+r_2(\widetilde H_{2*})+r_3(\widetilde H_{2*})\leq 2, \end{equation}
since the dual graph of $\overline B$ is a chain.

\subsection{Sufficient conditions.} Conversely, let us given a  rational Belyi pair $f:C\simeq \mathbb P^1\to\mathbb P_1$ of degree $n$ branched at $B_f\subset \{0,1,\infty\}$ whose cycle type of its monodromy is  $T(f)=\{ c_1,c_2,c_3\}$. Remind that the set of rational Belyi pairs is considered up to the actions of $PGL(2,\mathbb C)$ on $C$ and $\mathbb P^1$, and the ordered cycle type $T(f)=\{ c_1,c_2,c_3\}$ of $f$ depend on the choice of the base in $\pi_1(\mathbb P^1\setminus \{0,1,\infty\})$.
Therefore we can arrange the cycle type $T(f)$ so that the cycle type of $f_*(\overline{\gamma}_0)$ is $c_1$,  the cycle type of $f_*(\overline{\gamma}_{\infty})$ is $c_2$, and  the cycle type of $f_*(\overline{\gamma}_1)$ is $c_3$, where $\overline{\gamma}_0$,  $\overline{\gamma}_1$, $\overline{\gamma}_{\infty}$ are elements of $\pi_1(\mathbb P^1\setminus \{0,1,\infty\})$ represented by simple loops, resp., around $0$, $1$, and $\infty$ and such that $\overline{\gamma}_0 \overline{\gamma}_1 \overline{\gamma}_{\infty}=id$ in $\pi_1(\mathbb P^1\setminus \{0,1,\infty\})$.

\begin{df} We say that a  rational Belyi pair $f$ has the type: \newline
${\bf A}_{2k+1}$, $k\geq 1$, if  $(k_1-r_1)$  lengths $n_{1,j}$ in the cycle type $c_1=(n_{1,1},\dots, n_{1,k_1})$ are equal \newline $\phantom{aaaaa}$ to $k+1$, $r_1\leq 2$, and the last $r_1$ lengths are divisors of $k+1$;  \newline
${\bf A}_{2k}$, $k\geq 1$, if  $(k_1-r_1)$  lengths $n_{1,j}$ in the cycle type $c_1=(n_{1,1},\dots, n_{1,k_1})$ are equal to \newline $\phantom{aaaaa}$  $2k+1$, and the last $r_1$ lengths are divisors of $2k+1$,  \newline $\phantom{aaaaa}$ $(k_2-r_2)$  lengths $n_{2,j}$ in the cycle type $c_2=(n_{2,1},\dots, n_{2,k_2})$ are equal to $2$, the \newline $\phantom{aaaaa}$ last $r_2$ lengths are equal to $1$, and $r_1+r_2\leq 2$;  \newline
${\bf D}_{2k+2}$, $k\geq 2$, if  $(k_1-r_1)$  lengths $n_{1,j}$ in the cycle type $c_1=(n_{1,1},\dots, n_{1,k_1})$ are multi- \newline $\phantom{aaaaa}$  ple to $k$, $r_1\leq 2$;   \newline
${\bf D}_{2k+3}$, $k\geq 1$, if  $(k_1-r_1)$  lengths $n_{1,j}$ in the cycle type $c_1=(n_{1,1},\dots, n_{1,k_1})$ are multi- \newline $\phantom{aaaaa}$ ple to $2k+1$,   $(k_2-r_2)$ lengths $n_{2,j}$ in the cycle type $c_2=(n_{2,1},\dots, n_{2,k_2})$ are \newline $\phantom{aaaaa}$ equal to $2$, the last $r_2$ lengths are equal to $1$, and $r_1+r_2\leq 2$;
\newline
${\bf E}_{6}$ $\phantom{aa}$ if  $(k_1-r_1)$  lengths $n_{1,j}$ in the cycle type $c_1=(n_{1,1},\dots, n_{1,k_1})$ are equal to $4$, \newline $\phantom{aaaaa}$  and the last $r_1$ lengths are equal to $2$ or $1$,  $(k_2-r_2)$  lengths $n_{2,j}$ in the cycle \newline $\phantom{aaaaa}$ type $c_2=(n_{2,1},\dots, n_{2,k_2})$ are equal to $3$, the  last $r_2$ lengths are equal to $1$, and\newline $\phantom{aaaaa}$ $r_1+r_2\leq 2$; \newline
${\bf E}_7$ $\phantom{aa}$ if  $(k_1-r_1)$  lengths $n_{1,j}$ in the cycle type $c_1=(n_{1,1},\dots, n_{1,k_1})$ are even,  $(k_2-r_2)$ \newline $\phantom{aaaaa}$  lengths $n_{2,j}$ in the cycle type $c_2=(n_{2,1},\dots, n_{2,k_2})$ are equal to $3$, the last $r_2$ \newline $\phantom{aaaaa}$  lengths are equal to $1$, and $r_1+r_2\leq 2$;  \newline
${\bf E}_8$ $\phantom{aa}$ if  $(k_1-r_1)$  lengths $n_{1,j}$ in the cycle type $c_1=(n_{1,1},\dots, n_{1,k_1})$ are equal to $3$ and \newline $\phantom{aaaaa}$  the last $r_1$ lengths are equal to  $1$, $(k_2-r_2)$  lengths $n_{2,j}$ in the cycle type  $c_2=$  \newline $\phantom{aaaaa}$ $(n_{2,1},\dots, n_{2,k_2})$ are equal to $5$, the last $r_2$ lengths are equal to $1$, and $r_1+r_2\leq 2$.
\end{df}

\begin{rem} Note that a rational Belyi pair $f$ can have several $ADE$ types. For example if $f$ has ${\bf A}_{2k+1}$ type, $k\geq 1$, then it has also ${\bf D}_{2k+4}$ type. In addition, we will assume that any $f\in\mathcal{B}el$ has the ${\bf D}_4$ type.
\end{rem}

It follows from Lemmas \ref{A_n0} and \ref{3A_n}, the contractibility condition and inequality (\ref{nes}), that the {\it
necessary condition} for the branch curve $(B,o)$ of  a finite cover $F\in \mathcal R_T\subset \mathcal R\setminus (\mathcal R_{{\bf A}_0}\cup R_{{\bf A}_1})$ to belong to $\beta^{-1}(f)$, $\deg f=n$, is that $f$ {\it has type} $T$.

If this
necessary condition is met, then we can consider a monodromy homomorphism $H_{2*}:\pi_1^{loc}(B,o)\simeq \pi_1(\widetilde V\setminus \widetilde B)\to \mathbb S_n$ sending $\gamma_1$ to $f_*(\overline{\gamma}_0)$, $\gamma_2$ to $f_*(\overline{\gamma}_{\infty})$, $\gamma_3$ to $f_*(\overline{\gamma}_1)$, and $e$ to $id$. The homomorphism $H_{2*}$ defines finite coverings $H_2:(W,o_1)\to (V,o)$ and $\widetilde H_2:\widetilde W\to \widetilde V$, where $\sigma:\widetilde V\to (V,o)$ is the minimal resolution of the singular point of $(B,o)$. We can add to the maps $H_2$, $\widetilde H_2$, $\sigma$, and $\varsigma: \widetilde W\to (W,o_1)$ the bimeromorphic maps
$\varsigma_r:\overline W\to \widetilde W$,
$\varsigma_1:\overline W\to \overline W_m$ and
$\varsigma_m:\overline W_{m}\to (W,o_1)$.
As a result, we obtain the lower part of diagram $(*)$.
Denote by $C=\widetilde H_2^{-1}(E)$  the proper inverse image of the exceptional curve $E$ of the last $\sigma$-process. Obviously,  the restriction of $\widetilde H_2$ to $C$ coincides  with Belyi pair $f:C\to E$, $\deg f=\deg \widetilde H_2=n$.

As above, denote by  $m_0$ the order of the fundamental group $\overline{\pi}_1=\pi_1(W\setminus \{o_1\})\simeq\pi_1(\overline W_m\setminus \overline B^0)$, where $\overline B^0=\varsigma_m^{-1}(o_1)$. We obtain the following commutative diagram
\begin{equation}\label{diago}
\begin{CD} \overline U_{f,min_T} @>\overline{\varsigma}>> (U_{f,min_T},o_2) \\
 @V\overline H_f VV @VVH_f V \\
 \overline W_m @>\varsigma_m>> (W,o_1) \end{CD}
\end{equation}
in which $\overline H_f: \overline U_{f,min_T}\to\overline W_m$ and  $H_f: (U_{f,min_T},o_2)\to (W,o_1)$ are cyclic covers of degree $m_0$, $H_f:U_{f,min_T}\setminus \{o_2\}\to W\setminus \{ o_1\}$ is the universal cover, $\overline H_f$ is branched in $\overline B^0$,
$\overline U_{f,min_T}$ is a normal surface and $(U_{f,min_T},o_2)$ is a germ  of smooth surface, and $\overline{\varsigma}$ is a bimeromorphic holomorphic map.

Let $\gamma_{\overline C}\in\overline{\pi}_1$ be the element represented by a simple loop around $\overline C=\varsigma_1\circ \varsigma_r^{-1}(C)$.
We call the condition: {\it the dual graph of $\overline B^0$ is a chain and $\gamma_{\overline C}$ generates the group $\overline{\pi}_1$},
the {\it second necessary condition.}
If $f$ has type $T$ and the second necessary conditions is met for the germ $(B,o)$ having singularity type $T$, then $\overline H_f$  is branched in $\overline C$ with multiplicity $m_0$ and it is easy to see that the cover
$F_{f,min_T}:=H_2\circ H_f: (U_{f,min_T},o_2)\to (V,o)$
of degree $nm_0$ belongs to $\beta^{-1}(f)$. The cover $F_{f,min_T}\in \mathcal R_T 
$ will be called the {\it minimal cover in $\beta^{-1}(f)$ of  rational Belyi pair $f$ of type $T$.}

Denote by $\chi: \overline U\to \overline U_{f,min_T}$ the minimal resolution of singular points of $\overline U_{f,min_T}$. Then $\overline{\varsigma}\circ \chi:\overline U\to U_{f,min_T}$ is a composition of $\sigma$-processes with centers at nonsingular points. Note that $(\overline{\varsigma}\circ \chi)^{-1}(o_2)$ is a chain of the exceptional curves of $\overline{\varsigma}\circ \chi$. Denote by $\overline C_r=\chi^{-1}(\overline C)$ the proper inverse image of $\overline C$.

Consider the curve germ $B$ as a divisor in $(V,o)$ and let $F^*_{f,min_T}(B)=\sum r_jR_j$ be the inverse image of $B$, where $R_j$ are its irreducible components.
Denote by $S_1$ the set consisting of the pairs $(R_j,r_j)$ in which $R_j$ is a smooth germ for each $j$ and $(\overline{\varsigma}\circ \chi)^{-1}(R_j)$ is a chain consisting of the exceptional curves and the proper inverse image of $R_j$, and by $S_2$ the set consisting of the ordered pairs
$\{ (R_{j_1},r_{j_1}),(R_{j_2},r_{j_2})\}$, where $(R_{j_l},r_{j_l})\in S_1$ for $l=1,2$ and $R_{j_1}$, $R_{j_2}$ meet transversally at $o_2$.

The fundamental group $\pi_1(U_{f,min_T}\setminus R_j)\simeq \pi_1(\overline U\setminus (\overline{\varsigma}\circ \chi)^{-1}(R_j))$ is a free cyclic group generated by the element  $\gamma_j$ represented by a simple loop around $R_j$. Then
\begin{equation} \label{aj}\gamma_{\overline C}=\gamma_j^{a_j},\end{equation}
where $\gamma_{\overline C}$ is an element represented by a simple loop around $\overline C_r$ and $a_j$ can be computed using step by step Lemma \ref{sigma}.

Denote $M_{j}=\{ m\in \mathbb N\mid G.C.D.(m, a_j)=1\}$ for $(R_j,r_j)\in S_1$ and
$$M_{j_1,j_2}=\{ (m_1,m_2)\in\mathbb N^2\mid G.C.D.(m_1m_2, m_1a_{j_1}+m_2a_{j_2})=1\}$$
for $\{ (R_{j_1},r_{j_1}),(R_{j_2},r_{j_2})\}\in S_2$ and $a_j$ defined in (\ref{aj}) (note that if $(m_1,m_2)\in M_2$ then $G.C.D.(m_1,m_2)=1$).

For $\{ (R_{j_1},r_{j_1}),(R_{j_2},r_{j_2})\}\in S_2$ and $(R_{j_1},r_{j_1})\in S_1$ let us choose coordinates $(y_1,y_2)$ in $U_{f,min_T}$ such that $R_{j_l}$ is given by equation
$y_l=0$ for $l=1,2$ (if $(R_{j_1},r_{j_1})\in S_1$ then $R_{j_2}$ is any smooth curve germ meeting transversally with $R_{j_1}$)  and for each $(m_1,m_2)\in M_{j_1,j_2}$  (resp. for each $m_1\in M_{j_1}$ and $m_2=1$) consider the cyclic cover $\vartheta_{m_1,m_2}:(U,o')\to (U_{f,min_T},o_2)$ given by $x_1^{m_1}=y_1$, $x_2^{m_2}=y_2$.
It is easy to see that $F_{R_1,R_2,m_1,m_2}:=F_f\circ \vartheta_{m_1,m_2}:(U,o')\to (V,o)$ of degree $\deg F_{R_1,R_2,m_1,m_2}=nm_0m_1m_2$ also belongs to $\beta^{-1}(f)$.

Let $Aut(V,B,o)=\{ g\in Aut(V)\mid g(B)=B, g(o)=o\}$ be the automorphism group of the triple $(V,B,o)$ and $Gal(F_f)=\{ g\in Aut(U_{f,min_T},o_2)\mid F_f\circ g=F_f\}$ the automorphism group of $(U_{f,min_T},o_2)$ over $(V,o)$. The group $Aut(V,B,o)\times Gal(F_f)$ acts on the sets $S_1$ and $S_2$. Denote by $orb_j(f)$ the number of orbits of the action of the group $Aut(V,B,o)\times Gal(F_f)$  on $S_j$, $j=1,2$.

The results obtained above finally give the following
\begin{thm} \label{main2} Let $T$ be one of $ADE$ singularity types, a curve germ $(B,o)$  have the singularity of type $T$ at $o$, and $f\in \mathcal{B}el$, $\deg f=n$. Then the intersection $\mathcal R_T\cap \beta^{-1}(f)$ is nonempty if and only if $f$ has the type $T$ and $(B,o)$  satisfies the second necessary condition.

If $\mathcal R_T\cap \beta^{-1}(f)\neq \emptyset$ then $\mathcal R_T\cap \beta^{-1}(f)$ consists of the minimal cover $F_{f,min_T}:(U_{f,min_T},o_2)\to (V,o)$ of degree $nm_0$ and $(orb_1(f)+orb_2(f))$ infinite series of covers $F_{R_1,R_2,m_1,m_2}$, $\deg F_{R_1,R_2,m_1,m_2}=nm_0m_1m_2$.
\end{thm}

\section{Proof of Theorem \ref{B-D_4}} \label{thmD_4}

\begin{lem}\label{lem} Let a morphism $f:C=\mathbb P^1\to \mathbb P^1$, $\deg f=n>1$, be given by
$y_1=h_1(x_1,x_2)$, $y_2=h_2(x_2,x_2)$, where $h_1(x_1,x_2)$ and $h_2(x_2,x_2)$ are two coprime homogeneous in variables $x_1,x_2$ forms. Then the ramification divisor $R_f$ of $f$ is given by equation
\begin{equation} \label{J}
\displaystyle J_f(x_1,x_2):= \det  \left(\begin{array}{cc} \displaystyle
\frac{\partial h_1}{\partial x_1}, &
\displaystyle
\frac{\partial h_1}{\partial x_2} \\ \displaystyle
\frac{\partial h_2}{\partial x_1}, & \displaystyle
\frac{\partial h_2}{\partial x_2}\end{array}\right)= 0. \end{equation}
\end{lem}
\proof The projective line $C$ is covered by four neighbourhoods
$$U_{i,j}=\{ (x_1:x_2)\in C\mid x_i\neq  0,\,\, h_j(x_1,x_2)\neq 0\},\,\,\, 1\leq i,j\leq 2,$$
and $\widetilde x_{i}=\frac{x_{i}}{x_{\hat i}}$ is a coordinate in $U_{\hat i,j}$, where $\{i,\hat i\}=\{1,2\}$.
Similarly, $\mathbb P^1$ is covered by two affine lines $V_j=\{ (y_1:y_2)\in\mathbb P^1\mid y_j\neq 0\}$, $j=1,2$, and  $\widetilde y_{j}=\frac{y_{j}}{y_{\hat j}}$ is a coordinate in $V_{\hat j}$.  The morphism $f$ defines four rational functions $\widetilde y_j=f_{i,j}(\widetilde x_i)$ and, obviously, the restriction of $R_f$ to $U_{i,j}$ is the sum of critical points of $f_{i,j}$ counting with multiplicities. In particular, if $i=2$ and $j=2$ (the other cases are similar), then $R_f$ in $U_{2,2}$ is given by equation
$\frac{d\widetilde y_{1}}{d\widetilde x_1}=0$, where
$\widetilde y_1=f_{2,2}(\widetilde x_1)=\frac{h_1(\widetilde x_1,1)}{h_2(\widetilde x_1,1)}$. Therefore
\begin{equation} \label{pro} \frac{d\widetilde y_{1}}{d\widetilde x_1}=\frac{h_{1x_1}'(\widetilde x_1,1)h_2(\widetilde x_1,1)-h_1(\widetilde x_1,1)h_{2x_2}'(\widetilde x_1,1)}{h_2(\widetilde x_1,1)^2}.\end{equation}
It follows from Euler formula $nh(x_1,x_2)=x_1h'_{x_1}(x_1,x_2)+x_2h_{x_2}'(x_1,x_2)$ for homogeneous forms $h(x_1,x_2)$ of degree $n$ that
\begin{equation} \label{Eu} h_i(\widetilde x_1,1)=\frac{1}{n}[\widetilde x_1h'_{i\, x_1}(\widetilde x_1,1)+h_{i\, x_2}'(\widetilde x_1,1)] \end{equation}
and applying (\ref{Eu}), we obtain that  the numerator in right side of (\ref{pro}) coincides with $\frac{1}{n}J_f(\widetilde x_1,1)$. \qed

Applying direct calculations in non-homogeneous coordinates that define $\sigma$-processes with centers at points, it is easy to show that the following lemma occurs.
\begin{lem} \label{finite} Let $F:(U,o')\to (V,o)$ be a finite cover given by
\begin{equation}\label{generic} y_j=h_j(x_1,x_2)+\sum_{k=n_j+1}^{\infty}\sum_{m=0}^ka_{j,m}x_1^mx_2^{k-m},\quad j=1,2,\end{equation}
where $h_j(x_1,x_2)$ are homogeneous forms of degree $n_j$ in variables $x_1$ and $x_2$,
and let $\tau:\widetilde U\to U$ and $\sigma:\widetilde V\to V$ be $\sigma$-processes with centers at $o'$ and $o$, $\tau^{-1}(o')=\widetilde E$ and $\sigma^{-1}(o)=E$. Then \begin{itemize}
\item[$(i)$]  $\sigma^{-1}\circ F\circ \tau:\widetilde U\to\widetilde V$ is a holomorphic map if and only if $h_1(x_1,x_2)$ and $h_2(x_1,x_2)$ are coprime forms;
\item[$(ii)$] $\sigma^{-1}\circ F\circ \tau (\widetilde E)=E$ if and only if $n_1=n_2$ and the forms $h_1(x_1,x_2)$ and $h_2(x_1,x_2)$ are linear independent;
\item[$(iii)$] if $n_1=n_2:=n$ and  $h_1(x_1,x_2)$ and $h_2(x_1,x_2)$ are coprime forms, then
\begin{itemize} \item[$(iii)_1$] the restriction $f:\widetilde E\to E$ of $\sigma^{-1}\circ F\circ \tau$ to $\widetilde E$ is given by $y_1=h_1(x_1,x_2)$ and $y_2=h_2(x_1,x_2)$, \item[$(iii)_2$] $\deg F=n^2$ and $\deg f=n$. \end{itemize}
\end{itemize}
\end{lem}

To prove Theorem \ref{main2}, we use definitions and notations introduced in the previous sections.
We will assume that the branch locus $(B,o)$ of  a cover $F:(U,o') \to (V,o)$ is given by equation $uv(u-v)=0$. The cover $\widetilde H_2$ (see diagram ($*)$) is branched only in the disjoint union $B_1\sqcup B_2\sqcup B_3$ of three smooth curves, the proper inverse images of the irreducible branches of the curve $B$.  Therefore $\widetilde W$ is a smooth surface (that is, $\widetilde W=\overline W_m$) and the restriction of $\widetilde H_2$ to $C=\widetilde H_2^{-1}(E)$ is $\beta(F)$, $\deg \beta(F)=\deg \widetilde H_2:=n'$. Hence, $(C^2)_{\widetilde W}=-n'$ and $\varsigma: \widetilde W\to (W,o_1)$ is the minimal resolution of the singular point $o_1\in W$ of singularity type $A_{n',1}$ (that is $q=1$).
In addition, in diagram (\ref{diag-D_4}) the map $\overline{\varsigma}$ is the blowup of the point $\widetilde o$,
$\overline{\sigma}^{-1}(\widetilde o)=\widetilde{\theta}_{n',1}^{-1}(C)=\widetilde E$ and $(\widetilde E^2)_X=-1$; $\tau$ is the blowup of the point $o_1$,
$\tau^{-1}(o')=\widetilde{\vartheta}_{m_1,m_2}^{-1}(\widetilde E)=\overline E$ and $(\overline E^2)_{\widetilde U}=-1$.
In particular, all maps in (\ref{diag-D_4}) are holomorphic maps. Note also that
the restriction of
$\widetilde{\theta}_{n',1}$ to $\widetilde E$ is an isomorphism between  $\widetilde E$ and $C$. Therefore we can identify the restriction of
$F_{f,min}:=\widetilde H_2\circ \widetilde{\vartheta}_{n',1}$ to $\widetilde E$ with $\beta(F)$.

Let $F_{f,min}:(X,\widetilde o)\to (V,o)$ is given by equations (\ref{generic}). Then, by Lemma \ref{finite}, $f=\beta(F)$ is given by homogeneous forms  $y_1=h_1(x_1,x_2)$ and $y_2=h_2(x_2,x_2)$ of degree $n$. Therefore $n'=n$. In addition, according to Remark \ref{ident2}, we identify the monodromy homomorphism $\widetilde H_{2*}$ with monodromy homomorphism $f_*$.

If $F:(U,o')\to (V,o)$ is given by equations (\ref{b-d_4}), then $F_{f,min}$ is given by homogeneous forms $y_1=h_1(x_1,x_2)$ and $y_2=h_2(x_2,x_2)$. Therefore, first, we can  consider $F_{f,min}$ as the restriction to $\mathbb B_2\subset \mathbb C^2$ of the morphism $F_{f,min}:\mathbb C^2\to \mathbb C^2$ defined by the same functions. Second,
we can identify $C\simeq \mathbb P^1$ in the rational Belyi pair $f:C\to\mathbb P^1$ with the quotient space $\mathbb C^2/\{ (x_1,x_2)\sim (\lambda x_1,\lambda x_2)\,\,\text{for}\,\, \lambda\neq 0\}$ and the line $\mathbb P^1$ with $\mathbb C^2/\{ (y_1,y_2)\sim (\lambda y_1,\lambda y_2)\,\,\text{for}\,\, \lambda\neq 0\}$.
The ramification divisor $R_{F_{f,min}}$ of $F_{f,min}$ is given by equation (\ref{J}). Therefore $R_{F_{f,min}}$ is a sum of lines passing through the origin such that
$$R_{F_{f,min}}/\{ (x_1,x_2)\sim (\lambda x_1,\lambda x_2)\,\,\text{for}\,\, \lambda\neq 0\}=R_f$$
and it follows from Lemma \ref{lem} that the branch locus $B_{F_{f,min}}$ of $F_{f,min}$ is given by equation $uv(u-v)=0$ if $f\in \mathcal{B}el_3$. Therefore the singularity type of $B_{F_{f,min}}$ is ${\bf D}_4$. If $f\in \mathcal{B}el_2$ then  the branch locus $B_{F_{f,min}}$ of $F_{f,min}$ is given by equation $uv=0$ and therefore the singularity type of $B_{F_{f,min}}$ is ${\bf A}_1$. But, in both cases the branch locus $B_{F}$ of $F$ is given by equation $uv(u-v)=0$, since the branch curve of $\vartheta_{m_1,m_2}$ is contained in the union of two lines given by equations $x_1=0$, $x_2=0$ and $\{ f(p_1), f(p_2)\}\cup B_f=\{ 0,1,\infty\}$ for
$p_1=(0,1)$, $p_2=(1,0)\in C$.

Conversely, if $F\in \mathcal R_{{\bf D}_4}$ be given by equations (\ref{generic}), then it follows from above considerations that $f=\beta(F)$ is given by  homogeneous forms $y_1=h_1(x_1,x_2)$ and $y_2=h_2(x_2,x_2)$. Consider a cover $F'_{f,min}:(X',\widetilde o)\to (V,o)$ given by the same homogeneous forms $u=h_1(x_1,x_2)$ and $v=h_2(x_2,x_2)$ and consider  diagram (\ref{diag-D_4}) for $F'_{f,min}$ in which $\vartheta_{m_1,m_2}=\vartheta_{1,1}$ and we denote the germs of surfaces $W$, $\widetilde W$, $\widetilde X$ and the maps $H_2$, $\widetilde H_2$ and so on by the same letters with the addition of a stroke ($W'$, $\widetilde H_2'$ and so on). Then, by Lemmas \ref{lem} and \ref{finite}, $f'=\beta(F'_{f,min}):C'\to E$ is given by the same homogeneous forms $y_1=h_1(x_1,x_2)$ and $y_2=h_2(x_2,x_2)$.

According to Remark \ref{ident2}, the covers $\widetilde H_2$ and $\widetilde H'_2$ have the same monodromy homomorphism $f_*=f'_*$. Therefore, by Grauert - Remmert - Riemann - Stein Theorem, there are  biholomorphic isomorphisms  $\varphi:\widetilde W\to\widetilde W'$ and  $\psi:W\to W'$ such that $\widetilde H_2=\widetilde H_2'\circ \varphi$ and
$H_2=H_2'\circ \psi$, and we can identify $\widetilde W$ with $\widetilde W'$ and  $W$ with $W'$ with the help of these isomorphisms. Therefore there are  biholomorphic isomorphisms $\widetilde{\varphi}:\widetilde X\to \widetilde X'$ and $\widetilde{\psi}:X\to X'$ such that ${\varphi}\circ\widetilde{\theta}_{n,1}=\widetilde{\theta}'_{n,1}\circ\widetilde{\varphi}$ and ${\psi}\circ{\theta}_{n,1}={\theta}'_{n,1}\circ\widetilde{\psi}$, since $\widetilde{\theta}_{n,1}=\theta_{n,1}:\widetilde X\setminus \widetilde E=X\setminus\widetilde o\to \widetilde W\setminus C=W\setminus o_1$ and $\widetilde{\theta}'_{n,1}=\theta'_{n,1}:\widetilde X'\setminus \widetilde E'=X'\setminus\widetilde o'\to \widetilde W'\setminus C'=W'\setminus o'_1$ are the universal unramified covers. Hence, we can identify $\widetilde X$ with $\widetilde X'$ and  $X$ with $X'$. As a result, we obtain that the covers $F_{f,min}:X\to W$ and  $F'_{f,min}:X'\to W'$ are equivalent.

If $m_1$ or $m_2$ or both of them are grater than $1$ in $\vartheta_{m_1,m_2}:(U,o')\to (X,\widetilde o)$, then we can obtain the cyclic cover $\vartheta'_{m_1,m_2}:(U',o'')\to(X',\widetilde o')$ branched in the images of the branch curves of $\vartheta_{m_1,m_2}$ under the holomorphic isomorphism $\widetilde{\psi}$. Obviously, the covers $\vartheta_{m_1,m_2}$ and $\vartheta'_{m_1,m_2}$ are equivalent. Therefore the compositions $F=F_{f,min}\circ\vartheta_{m_1,m_2}$ and $F'=F'_{f',min}\circ\vartheta'_{m_1,m_2}$ of equivalent covers are equivalent and it is easy to see that there is a coordinate change in $(U',o'')$ such that the cover $F'$ is given by equations  (\ref{b-d_4}).

\section{Proof of Theorem \ref{beta2}} \label{Bel2}

\subsection{} Let $f=\beta (F)\in \mathcal{B}et_2$, $\deg f=n$, for $F$ branched in $(B,o)$ having one of $ADE$ singularity types. Without loss of generality we can assume that $f$ is given by $y=x^n$ in non-homogeneous coordinates and its branch locus is $B_f=\{ 0,\infty\}\subset \{ 0,1, \infty\}$. Then, in general case, $\widetilde H_2:\widetilde W\to \widetilde V$ is a cyclic cover  branched in two of three trails of $\widetilde B$,
$\deg \widetilde H_2=n$, the cover $\widetilde H_1=\widetilde{\theta}_{n,q}\circ\widetilde{\vartheta}_{m_1,m_2}:\widetilde U\to \widetilde W$ is also a cyclic cover, and $\widetilde{\vartheta}_{m_1,m_2}$ must be branched at least in one of the irreducible components of the inverse image of the third trail.

\subsection{Cases $\mathcal R_{{\bf A}_{0}}$ and $\mathcal R_{{\bf A}_{1}}$.} 
It is easy to show that if $F\in (\mathcal R_{{\bf A}_0}\cup\mathcal R_{{\bf A}_1})\cap\beta^{-1}(\mathcal{B}el)$, $\deg \beta(F)=n$, then $F$ is equivalent one of the following germs of covers of degree $n^2m_1m_2$ given by
$u=z^{nm_1}$, $v=w^{nm_2}$ for some $m_1\geq m_2\geq 1$, $G.C.D.(m_1,m_2)=1$ {\rm (}if $F\in \mathcal R_{{\bf A}_0}$, then $n=m_2=1${\rm )}.

\subsection{Case $\mathcal R_{{\bf A}_{2k+1}}$, $k\geq 1$.} \label{case} Without loss of generality, we can assume that $(B,o)$ is given by equation $u(u-v^{k+1})=0$ and  $u=0$ is an equation of the irreducible component $B_1$ of $B$.
The graph $\Gamma(\widetilde B)$ of the curve germ $(B,o)$ of singularity type ${\bf A}_{2k+1}$, $k\geq 1$, is depicted on Fig. 1 (in this case $E=E_{k+3}$). Let us number the trails of $\widetilde B$ as follows: $\widetilde B_1=B_1$, $\widetilde B_2=B_2$, and $\widetilde B_3=E_{3}\cup\dots\cup E_{k+2}$. Then $\widetilde H_2$ is branched either in $B_1\cup B_2$ or in $B_l\cup (\bigcup_{j=3}^{k+2}E_j)$, where $l=1$ or $2$.

Let us show that the first case is impossible. Indeed, if $\widetilde H_2$ is branched in $B_1\cup B_2$, then $n$ must be equal to $2$, since
$\widetilde H_2^{-1}(\bigcup_{j=3}^{k+3}E_j)$ must be a chain of rational curves satisfying the second necessary condition. But, if $n=2$ then $\widetilde H_1$ must be branched only  in $\widetilde H_2^{-1}(\bigcup_{j=3}^{k+3}E_j)=\bigcup_{j=1}^{2k+1}\widetilde E_j$, where $\widetilde E_{k+1}=\widetilde H_2^{-1}(E_{k+3})=C$. In this case, $\widetilde W$ is a smooth surface and $(\widetilde E_j^2)_{\widetilde W}=-2$ for all $j$. Therefore, by Lemma \ref{A_n0}, the group $\pi_1(N_T\setminus (\bigcup_{j=1}^{2k+1}\widetilde E_j))\simeq \mu_{2k+2}$, where $N_T$ is a small tubular neighbourhood of $\bigcup_{j=1}^{2k+1}\widetilde E_j$. It follows from Theorem \ref{Mum} that the group $\pi_1(N_T\setminus (\bigcup_{j=1}^{2k+1}\widetilde E_j))$ is generated by an element $\gamma_{\widetilde E_1}$ represented by a simple loop around $\widetilde E_1$ and $\gamma_{C}=\gamma_{\widetilde E_1}^{k+1}$, where $\gamma_C$ is an element represented by a simple loop around $\widetilde E_{k+1}$. Therefore in this case the second necessary condition is not satisfied, since $\gamma_{C}=\gamma_{\widetilde E_1}^{k+1}$ does not generate the group $\mu_{2k+2}$.

In the second case without loss of generality, we can assume that $\widetilde H_2$ is branched in $\widetilde B_1= B_1$ and $\widetilde B_3=\bigcup_{j=3}^{k+2}E_j$, and in addition, assume that $n=\deg \widetilde H_2$ is a divisor of $k+1$, $k+1=nm_0$, since $\widetilde{\pi}^0_3\simeq \mu_{k+1}$.

Consider the surface $\widetilde W_m$ and the curve $\overline B\subset \overline W_m$ (see diagram $(*)$). The inverse image
$(H_2\circ \varsigma_m)^{-1}(B_2)=\bigsqcup_{j=1}^n\overline B_{2,j}$ is the disjoint union of $n$ curve germs each of which intersects transversally with $\overline C$ and $(H_2\circ \varsigma_m)^{-1}(B_1)=\overline B_{1,1}$ is an irreducible curve germ.
By Lemma \ref{delta1}, there are three possibilities for the curve $\overline B^0\subset \overline W_m$ (see notations introduced  in Subsection \ref{neces}). The first one (when $m_0>1$) is
$\overline B^0=\overline C\cup (\bigcup_{j=1}^{m_0-1}\overline E_j)$, the curve $\overline B_1$ is one of the irreducible components of $(H_2\circ \varsigma_m)^{-1}(B_2)$ (say $\overline B_{2,1}$), and   $\overline B_2=\emptyset$. In this case we have $(\overline E_j^2)_{\overline W_m}=-2$ and $(\overline C^2)_{\overline W_m}=-1$. In the second case (when $m_0=1$)
$\overline B_1=\overline B_{1,1}$, $\overline B_2$ is one of the irreducible components of $(H_2\circ \varsigma_m)^{-1}(B_2)$, and
$\overline B^0=\overline C$, $(\overline C^2)_{\overline W_m}=-1$.
In the third case (when $m_0=1$)
$\overline B_1$ and $\overline B_2$ are two of the irreducible components of $(H_2\circ \varsigma_m)^{-1}(B_2)$, and
$\overline B^0=\overline C$, $(\overline C^2)_{\overline W_m}=-1$.

In all of these cases, $(W,o_1)$ is a germ of smooth surface and $H_2$ is given by $u=y_1^n$ and $v=y_2$, where $y_1$, $y_2$ are coordinates in $(W,o_1)$. We have
$H_2^{-1}(B_2)=\bigcup_{j=1}^nB_{2,j}$, where $B_{2,j}$ are given by equations $y_1-\omega_jy_2^{m_0}=0$, $\omega_j=exp(2\pi j i/n)$. May be after scalar coordinate changers in $(V,o)$ and $(W,o_1)$, we can assume that one of irreducible components of $H_2^{-1}(B_2)$ included in the branch curve of $H_1$ is given by equation $y_1-y_2^{m_0}=0$.

In the first case, put $x_1=y_1-y_2^{m_0}$ and $x_2=y_2$. Then it is easy to see that $H_1:(U,o')\to (W,o_1)$ is given by functions $z^m=x_1$ and $w=x_2$, where $z$, $w$ are coordinates in $(U,o')$ and $m>1$. In view of Lemma \ref{sigma}, the second necessary condition entails the equality $G.C.D.(m_0,m)=1$. Therefore $F:(U,o')\to (V,o)$ is given by
$$ u=(z^m+w^{m_0})^n,\quad v=w,$$
where $n,m,m_0>1$ and $G.C.D.(m_0,m)=1$.

In the second case, put $x_1=y_1$ and $x_2=y_2-y_1$. Then it is easy to see that $H_1:(U,o')\to (W,o_1)$ is given by functions $z^{m_1}=x_1$ and $w^{m_2}=x_2$, where $z$, $w$ are coordinates in $(U,o')$ and $m_1\geq 1$, $m_2>1$,  $G.C.D.(m_1,m_2)=1$. Therefore $F:(U,o')\to (V,o)$ is given by
$$ u=z^{nm_1},\quad v=z^{m_1}+w^{m_2},$$
where $m_1\geq 1$, $n,m_2>1$, and $G.C.D.(m_1,m_2)=1$.

In the third case, the cover $H_1:(U,o')\to (W,o_1)$ is branched in two curves given by $y_1-y_2=0$ and $y_1-\omega_jy_2=0$ for some $j$, $1\leq j\leq n-1$. We
put $x_1=(\omega_j-1)^{-1}(y_1-y_2)$ and $x_2=(\omega_j-1)^{-1}(y_1-\omega_jy_2)$. Then $H_1:(W,o_1)\to (V,o)$ is given by functions $z^{m_1}=x_1$ and $w^{m_2}=x_2$, where $z$, $w$ are coordinates in $(U,o')$ and $m_1$, $m_2>1$,  $G.C.D.(m_1,m_2)=1$. Therefore $F:(U,o')\to (V,o)$ is given by
$$ u=(\omega_jz^{m_1}-w^{m_2})^n,\quad v=z^{m_1}-w^{m_2},$$ where $n, m_1,m_2>1$, $G.C.D.(m_1,m_2)=1$, and $\omega_j={\bf e}^{2\pi j i/n}$, $1\geq j\geq n-1$.

\subsection{Case $\mathcal R_{{\bf A}_{2k}}$, $k\geq 1$.} \label{caseA_2k}
The graph $\Gamma(\widetilde B)$ of the curve germ $(B,o)$ of singularity type ${\bf A}_{2k}$, $k\geq 1$, is depicted on Fig. 2 (in this case $E=E_{k+2}$). Let us show that in this case $(\bigcup_{k=1}^{\infty}\mathcal R_{{\bf A}_{2k}})\cap \beta^{-1}(\mathcal{B}el_2)=\emptyset$. Indeed, assume that there is $F\in \mathcal R_{{\bf A}_{2k}}\cap \beta^{-1}(\mathcal{B}el_2)$ for some $k\geq 1$. Then (see diagram $(*)$) $\widetilde H_2$ can be branched either in $\widetilde B_1\cup \widetilde B_2$, or in
$\widetilde B_1\cup \widetilde B_3$, or in $\widetilde B_2\cup\widetilde B_3$, where  $\widetilde B_1=E_{k+3}$, $\widetilde B_2=E_{2}\cup\dots\cup E_{k+1}$, and $\widetilde B_3=B_1$.

By Lemmas \ref{A_n0} and \ref{3A_n}, $\widetilde{\pi}^0_1\simeq \mu_2$ and $\widetilde{\pi}_2\simeq \mu_{2k+1}$. Therefore $\widetilde H_2$ can not be branched in
$\widetilde B_1\cup \widetilde B_2$, since $G.C.D.(2, 2k+1)=1$.

Assume that $\widetilde H_2$ is branched in $\widetilde B_2\cup \widetilde B_3$. Then $\deg \widetilde H_2$ is a divisor of $2k+1$ and hence the dual graph of
$\widetilde H_2^{-1}(E_{k+2}\cup E_{k+3})$ is not a tree, that is, the second necessary condition is not satisfied.

Assume that $\widetilde H_2$ is branched in $\widetilde B_1\cup \widetilde B_3$. Then $\deg \widetilde H_2=2$ and the dual graph of $\widetilde H_2^{-1}(\widetilde B_{2}\cup E_{k+2})=\bigcup_{j=1}^{2k+1}\widetilde E_j$ (here $\widetilde E_{k+1}=\widetilde H_2^{-1}(E)$=C) is a chain with weights $[\underbrace{2,\dots,2}_{k-1},3,2,3,\underbrace{2,\dots,2}_{k-1}]$. Therefore $\overline B^0=\bigcup_{j=1}^{2k+1}\overline E_j\subset \overline W_m$ is a chain with the weights $[\underbrace{2,\dots,2}_{k-1},3,1,3,\underbrace{2,\dots,2}_{k-1}]$, since $(\widetilde H_2^{-1}(E_{k+3}),\widetilde H_2^{-1}(E_{k+3}))_{\widetilde W}=-1$. It follows from Theorem \ref{Mum} (or see the proof of Lemma \ref{delta2}) that $\gamma_{\overline C}=1$ in $\pi_1(\overline W_m\setminus \overline B^0)$, that is,  the second necessary condition is not satisfied in this case.

\subsection{Case $\mathcal R_{{\bf D}_{2k+3}}$, $k\geq 1$.}
The graph $\Gamma(\widetilde B)$ of the curve germ $(B,o)$ of singularity type ${\bf D}_{2k+3}$, $k\geq 1$, is depicted on Fig. 4 (in this case $E=E_{k+3}$). We can assume that $(B,o)$ is  given by equation $u(v^2-u^{2k+1})=0$.  The cyclic cover $\widetilde H_2:\widetilde W\to\widetilde V$  (see diagram $(*)$) can be branched either in $\widetilde B_1\cup \widetilde B_2$, or in
$\widetilde B_1\cup \widetilde B_3$, or in $\widetilde B_2\cup\widetilde B_3$, where  $\widetilde B_1=B_1\cup E_{3}\cup\dots\cup E_{k+2}$, $\widetilde B_2=E_{k+4}$, and $\widetilde B_3=B_2$, but using the same arguments as in Subsection \ref{caseA_2k}, it is easy to show that $\widetilde H_2$ can be branched neither in  $\widetilde B_1\cup\widetilde B_3$ nor in $\widetilde B_2\cup\widetilde B_3$.

By Remark \ref{gener} and Lemmas \ref{A_n0} and \ref{3A_n}, $\widetilde{\pi}_1$ is the infinite cyclic group generated by $e_{k+2}$ and $\widetilde{\pi}^0_1\simeq \mu_{2k+1}$,  $\widetilde{\pi}^0_2\simeq \mu_{2}$.
Therefore if $\widetilde H_2$ is branched in
$\widetilde B_1\cup \widetilde B_2$, then $\deg \widetilde H_2=2$ and $\widetilde H_{2*}(e_{k+2})$ generates the monodromy group $G_{\widetilde H_2}\simeq \mu_2$. Applying the representation of the group $\widetilde{\pi}_1$  obtained with the help of Theorem \ref{Mum}, it is not difficult to show that $\widetilde W$ has singular points of type $A_1$ over the intersection points of the irreducible components of $\widetilde B$ and $\overline B^0\subset \overline W_m$ is a chain of rational curves, $\overline B^0=\overline C\cup (\bigcup_{j=1}^{2k}\overline E_j$, the weights of its dual graph are $[\underbrace{2,\dots,2}_{2k},1]$. Therefore $(W,o_1)$ is a germ of smooth surface and $H_2:(W,o_1)\to (V,o)$ is given by functions $y_1^2=u$ and $y_2=v$, where $y_1$, $y_2$ are coordinates in $(W,o_1)$. The inverse image $H_2^{-1}(B_2)=B_{2,1}\cup B_{2,2}$ is the union of two branches given by equations $y_2-y_1^{2k+1}=0$ and $y_2+y_1^{2k+1}=0$. The cyclic cover $H_2:((U,o')\to (W,o_1)$ is branched  in $H_2^{-1}(B_1)$ given by equation $y_1=0$ with multiplicity $m_1$ and in one of irreducible components of $H_2^{-1}(B_2)$ with multiplicity $m_2$. Without loss of generality, we can assume that $H_1$ is branched in $B_{2,1}$. Put $x_1=y_1$ and $x_2=y_2-y_1^{2k+1}$. Then $H_1$ is given by functions $z^{m_1}=x_1$ and $w^{m_2}=x_2$. Therefore $F:(U,o')\to (V,o)$ is given by functions
$$ u=z^{2m_1}, \,\, v=z^{m_1(2k+1)}+w^{m_2},$$
where $k$, $m_1\geq 1$, $m_2>1$ and $G.C.D.(m_1,m_2)=1$.
Note that $\varsigma_m:\overline W_m\to (W,o_1)$ is a composition of $2k+1$ $\sigma$-processes blowing up the point $o_1$ and points lying in the proper inverse images of $B_{2,2}$ and such that $\overline C$ is the exceptional curve of the last blowup. Therefore it follows from the second necessary condition that $G.C.D.(m_2,2k+1)=1$.

\subsection{Cases $\mathcal R_{{\bf E}_{6}}$ and $\mathcal R_{{\bf E}_{8}}$.}
The graphs $\Gamma(\widetilde B)$ of the curve germs $(B,o)$ of singularity types ${\bf E}_{6}$ and ${\bf E}_8$ are depicted on Fig. 6 and 8. We will show that
$\mathcal R_{{\bf E}_{6}}\cap \beta^{-1}(\mathcal{B}el_2)=\emptyset$ (the proof that $\mathcal R_{{\bf E}_{8}}\cap \beta^{-1}(\mathcal{B}el_2)=\emptyset$ is similar and therefore it will be omitted). Assume that there is $F\in \mathcal R_{{\bf E}_{6}}\cap \beta^{-1}(\mathcal{B}el_2)$. Then (see diagram $(*)$) $\widetilde H_2$ can be branched either in $\widetilde B_1\cup \widetilde B_2$, or in
$\widetilde B_1\cup \widetilde B_3$, or in $\widetilde B_2\cup\widetilde B_3$, where  $\widetilde B_1=E_{2}$, $\widetilde B_2=E_{4}\cup E_{5}$, and $\widetilde B_3=B_1$.

It follows from Theorem  \ref{Mum} that $\widetilde{\pi}^0_1\simeq \mu_4$ and $\widetilde{\pi}_2\simeq \mu_{3}$. Therefore $\widetilde H_2$ can not be branched in
$\widetilde B_1\cup \widetilde B_2$, since $G.C.D.(4, 3)=1$.

If $\widetilde H_2$ is branched in $\widetilde B_1\cup \widetilde B_3$ or in $\widetilde B_2\cup \widetilde B_3$, then it is easy to see that the dual graph of
$\overline B^0\subset \overline W_m$  is not a tree, that is, the second necessary condition is not satisfied.

\subsection{Case $\mathcal R_{{\bf E}_{7}}$.}
The graph $\Gamma(\widetilde B)$ of the curve germ $(B,o)$ of singularity type ${\bf E}_{7}$ given by $u(v^2-u^3)=0$ is depicted on Fig. 6 (in this case $E=E_{4}$). Let us number the trails of $\widetilde B$ as follows: $\widetilde B_1=B_1\cup E_5$, where $B_1$ is given by equation $u=0$, $\widetilde B_2=E_3$, and $\widetilde B_3=B_2$. Similar to the cases considered above, it is easy to see that $\widetilde H_2$ must be branched in $B_1\cup B_2$, $\deg \widetilde H_2=3$, $(W,o_1)$ is a germ of smooth surface and $H_2$ is given by $u=y_1^3$ and $v=y_2$. Then $H_1^{-1}(B_2)=B_{2,1}\cup B_{2,2}\cup B_{2,3}$, where $B_{2,j}\subset W$ are given in coordinates $y_1$, $y_2$ in $(W,o_1)$ by equations $y_2-\omega_jy_1^2=0$, $\omega_j={\bf e}^{2\pi ji/3}$, $j=1,2,3$. The cover  $H_1:(U,o')\to (W,o_1)$ is branched in one of irreducible components of $H_2^{-1}(B_2)$ with multiplicity $m_2>1$ and, possibly, in $H_2^{-1}(B_1)$. As above, without loss of generality, we can assume that  $H_1:(U,o')\to (W,o_1)$ is branched in $B_{2,3}$.
Put $x_1=y_1$ and $x_2=y_2-y_1^2$. Then  $H_1$ 
is given by functions $z^{m_1}=x_1$ and $w^{m_2}=x_2$, where $z$, $w$ are coordinates in $(U,o')$ and  $G.C.D.(m_1,m_2)=1$. Therefore $F:(U,o')\to (V,o)$ is given by
$$ u=z^{3m_1},\quad v=z^{2m_1}+w^{m_2},$$
where $m_1\geq 1$, $m_2>1$, and $G.C.D.(m_1,m_2)=1$.

\subsection{Case $\mathcal R_{{\bf D}_{4}}$.} It easily follows from the proof of Theorem \ref{thmD_4} that a germ of cover $F\in \mathcal R_{{\bf D}_4}\cap \beta^{-1}(\mathcal{B}el_2)$ is equivalent to a cover given by one of the following pairs of functions:
\begin{equation} u=z^{m_1n},\,\, v=(z^{m_1}+w^{m_2})^n \end{equation}
or
\begin{equation} u=(z^{m_1}-w^{m_2})^n,\,\, v=(z^{m_1}-\omega_jw^{m_2})^n, \end{equation}
where $n\geq 2$, $m_1$, $m_2\geq 1$ and $G.C.D.(m_1,m_2)=1$, $\omega_j=exp(2\pi ji/n)$, $1\leq j\leq n-1$.

\subsection{Case $\mathcal R_{{\bf D}_{2k+2}}$, $k\geq 2$.} Without loss of generality, we can assume that $(B,o)$ is given by equation $uv(v-u^{k})=0$, where  $u=0$ is an equation of the irreducible component $B_1$ of $B$ and $v=0$ is an equation of the irreducible component $B_2$.
The graph $\Gamma(\widetilde B)$ of the curve germ $(B,o)$ of singularity type ${\bf D}_{2k+2}$, $k\geq 2$, is depicted on Fig. 3 (in this case $E=E_{k+3}$). Let us number the trails of $\widetilde B$ as follows: $\widetilde B_1=B_1\cup E_4\cup \dots\cup E_{k+2}$, $\widetilde B_2=B_2$, and $\widetilde B_3=B_{3}$. The cyclic cover $\widetilde H_2:\widetilde W\to \widetilde V$
is branched either in $\widetilde B_1\cup \widetilde B_2$ or in $\widetilde B_2\cup\widetilde B_3$ (the case when $\widetilde H_2$ is branched in $\widetilde B_1\cup \widetilde B_3$ is the same as the case when $\widetilde H_2$ is branched in $\widetilde B_1\cup \widetilde B_2$, since we can make a coordinate change in $(V,o)$).

As in Subsection \ref{case}, it is easy to show that the case when $\widetilde H_2$ is branched in $\widetilde B_2\cup \widetilde B_3$ is impossible.

Consider the case when the cyclic cover $\widetilde H_2:\widetilde W\to \widetilde V$ (see diagram $(*)$)
is branched in $\widetilde B_1\cup \widetilde B_2$, $\deg \widetilde H_2=n$. Let $n=n_1k_1$ and $k=k_1k_2$, where $k_1=G.C.D.(n,k)$ and $G.C.D.(n,k_2)=1$. The group $\widetilde{\pi}_1$ is generated by element $\gamma_{E_{k+2}}$ represented by a simple loop around the curve $E_{k+2}$. It follows from Theorem \ref{Mum} that   $\gamma_{B_1}=\gamma^k_{E_{k+2}}$, where $\gamma_{B_1}$ is an element in $\widetilde{\pi}_1$ represented by a simple loop around the germ $B_{1}$. The monodromy group $G_{\widetilde H_2}\simeq \mu_n\subset\mathbb S_n$ is generated by element $g=\widetilde H_{2*}(\gamma_{E_{k+2}})$. The element $\widetilde H_{2*}(\gamma_{B_2})=g^{-1}$ is also a generator of  the group $\mu_n$, where  $\gamma_{B_2}$ is an element in $\widetilde{\pi}_1$ represented by a simple loop around the germ $B_{2}$, and  $\widetilde H_{2*}(\gamma_{B_1})=g^{k}$ is an element of order $n_1$. Therefore $\widetilde H_2$ is branched in $B_1$ with multiplicity $n_1$ and in $B_2$ with multiplicity $n$. As a result, we obtain that $H_2$ also is branched only in $B_1$ with multiplicity $n_1$ and in $B_2$ with multiplicity $n$.

In diagram (\ref{diag-D_4}), the cover $\theta_{n',q}:(X,\widetilde o)\to (W,o_1)$ is branched only at the point $o_1$. Therefore the map $F_{f,min}=H_2\circ \theta_{n',q}:(X,\widetilde o)\to (V,o)$ is given in some coordinates $(x_1,x_2)$ in $(X,\widetilde o)$ by functions
$$ u=x_1^{n_1},\quad v=x_2^n,$$
since $(X,\widetilde o)$ is a germ of smooth surface and $H_2\circ \theta_{n',q}$ is branched in the divisor with normal crossing $B_1\cup B_2$. The inverse image $F^{-1}_{f,min}(B_3)=\cup_{j=1}^nB_{3,j}$ is the union of $n$ smooth curves given by equation $x_2^n-x_1^{n_1k}=\prod_{j=1}^n(x_2-\omega_jx_1^{k_2})=0$, where $\omega_j=exp(2\pi ji/n)$, $j=1,\dots, n$.

The map $\vartheta_{m_1,m_2}:(U,o')\to (X,\widetilde o)$ is branched in at most two irreducible curves, one of which belongs to $F^{-1}_{f,min}(B_3)$. Without loss of generality, we can assume that it is given by equation $y_2:=x_2-x_1^{k_2}=0$. If the branch locus of $\vartheta_{m_1,m_2}$ consists of two irreducible components, then the other one is an irreducible component of the inverse image either of $B_1$, or $B_2$, or $B_3$. Therefore we have three possibilities: the second irreducible component is given by equation either $y_1:=x_1=0$, or $y_1:=x_2=0$ (if $k_2=1$), or $y_1:=x_2-\omega_jx_1^{k_2}=0$ for some $j=1,\dots, n-1$ (if $k_2=1$), and $\vartheta_{m_1,m_2}$ is given by functions $y_1=z^{m_1}$ and $y_2=w^{m_2}$. Applying Lemma \ref{sigma}, the second necessary condition is equivalent to the condition: $G.C.D.(m_1,m_2)=G.C.D.(m_2,k_2)=1$ in the first case and   $G.C.D.(m_1,m_2)=1$ in the second and third cases. As a result, we obtain that $F:(U,o')\to (V,o)$ is equivalent to one of the following covers given by functions:
$$ u=z^{n_1m_1},\quad v=(z^{m_1k_2}+w^{m_2})^{n_1k_1}
$$
in the first case;
$$ u=(z^{m_1}-w^{m_2})^{n_1},\quad v=z^{m_1n_1k_1}
$$
in the second case;
$$ u=(z^{m_1}-w^{m_2})^{n_1},\quad v=(z^{m_1}-\omega_jw^{m_2})^{n_1k_1}
$$
 in the third case, where $k_1k_2\geq 2$, $n_1k_1\geq 2$, $m_1,m_2\geq 1$, and $G.C.D.(m_1,m_2)=G.C.D.(nm_2,k_2)=1$, $\omega_j=exp(2\pi ji/n_1k_1)$, $j=1,\dots, n_1k_1-1$.

\end{document}